\documentclass[11pt]{amsart}

\usepackage{amssymb,amsthm,amsmath,eucal,mathrsfs}

\usepackage{enumitem}
\usepackage{hyperref}

\advance\textheight by \topskip

\def\N{\mathbb{N}}
\def\R{\mathbb{R}}

\def\C{\mathbb{C}}

\def\E{\mathcal{E}}
\def\cL{\mathcal{L}}

\def\PP{\mathbb{P}}
\def\EE{\mathbb{E}}

\renewcommand{\check}{\breve}

\def\dv{{\mathrm{div}\,}}
\def\rot{{\mathrm{rot}\,}}

\def\loc{{\mathrm{loc}}}

\def\dive{\mbox{\rm div}\,}

\def\arg{\mbox{\rm arg}\,}

\def\sm{\setminus}

\def\bel{\begin{equation}}
\def\eel{\end{equation}}

\def\beq{\begin{eqnarray*}}
\def\eeq{\end{eqnarray*}}

\newcommand{\sk}[2]{\langle #1,#2\rangle}

\newcommand{\n}[1]{\|#1\|}
\newcommand{\nn}[2]{\|#1\|_{#2}}

\def\Tr{\mathrm{Tr}}
\def\RTr{R_{\mathrm{Tr}}}

\def\eps{\varepsilon}
\def\ph{\varphi}
\def\Om{\Omega}
\def\om{\omega}
\def\wt{\widetilde}
\def\wh{\widehat}
\def\la{\lambda}

\def\al{\alpha}
\def\ga{\gamma}

\def\si{\sigma}
\def\Si{\Sigma}
\def\del{\delta}

\def\ov{\overline}

\def\vrh{\varrho}
\def\a{\mathfrak{a}}
\def\Rep{\mathrm{Re}\,}
\def\Imp{\mathrm{Im}\,}

\def\emb{\hookrightarrow}

\newtheorem{prop}{Proposition}[section]
\newtheorem{lemma}[prop]{Lemma}
\newtheorem{remark}[prop]{Remark}
\newtheorem{corollary}[prop]{Corollary}
\newtheorem{definition}[prop]{Definition}
\newtheorem{theorem}[prop]{Theorem}

\newtheorem{assumption}[prop]{Assumption}

\title{$H^\infty$-calculus for the Stokes operator with Hodge, Navier, and Robin boundary conditions on unbounded domains}

\author{Peer Christian Kunstmann}
\address{\it Karlsruhe Institute of Technology (KIT),
Institute for Analysis\\
Englerstr. 2, D -- 76128 Karlsruhe, Germany\\
e-mail: peer.kunstmann@kit.edu}

\date{}

\begin{document}

\begin{abstract}
We study the Stokes operator with Hodge, Navier, and Robin boundary conditions on domains $\Om\subseteq\R^d$ that are uniformly $C^{2,1}$. Starting with the Hodge Laplacian we etablish a bounded H\"ormander functional calculus for the Stokes operator with Hodge boundary conditions. This entails a H\"ormander functional calculus and boundedness of the $H^\infty$-calculus in spaces of soleniodal vector fields for the Stokes operator with Hodge boundary conditions. We then establish boundedness of the $H^\infty$-calculus for Stokes operators with Navier type conditions via Robin type perturbations of Hodge boundary conditions. This implies maximal $L^p$-regularity for these operators and results on fractional domain spaces.
Our results cover certain non-Helmholtz domains.


\end{abstract}

\maketitle

\section{Introduction}\label{sec:intro}

Boundary condtions of Navier type play a vital role in the mathematical investigation of problems in fluid mechanics. They are used to model various slip type condtions on a fixed wall. In this paper we study Stokes operators on unbounded uniform $C^{2,1}$-domains under Hodge (also called perfect slip) conditions and boundary conditions of Navier type. 

It is well-known that, for general $C^{2,1}$-domains $\Om\subseteq\R^d$, the Helmholtz decomposition of $L^q(\Om)^d$ into the solenoidal space $L^q_\si(\Om)$ and the gradient space $G^q(\Om)$ may fail for certain $q\in(1,\infty)$, see \cite{MaBo}. 
As a way out,   
the spaces $\wt{L}^q_\si(\Om)$ have been introduced by Farwig, Kozono, and Sohr (see, e.g., \cite{FKS:ArchMath}). On the other hand, there has been an interest in recent years in the study of Stokes and Navier-Stokes equations also in certain unbounded non-Helmholtz domains. In particular, the results on Stokes operators with Navier type boundary conditions by Hobus and Saal in \cite{hobus-saal} cover situations in which the Helmholtz decomposition of $L^q(\Om)^d$ fails.    

Following \cite{hobus-saal} we treat Navier type boundary conditions as a special case of Robin type perturbations of Hodge boundary conditions. Under certain assumptions on the domain $\Om$ and $q\in(1,\infty)$, it has been shown in \cite{hobus-saal} that Stokes operators with Hodge and Navier type boundary conditions generate analytic semigroups. For the spaces $\wt{L}^q_\si(\Om)$, $1<q<\infty$, it has been shown by Farwig and Rosteck in \cite{FarRost-Res}, \cite{FarRost-MR}, that Stokes operators with Navier type boundary conditions generate analytic semigroups and have the property of maximal $L^p$-regularity, $1<p<\infty$. 
In this paper we substantially extend these results by establishing a bounded $H^\infty$-calculus and maximal $L^p$-regularity 
for Stokes operators with Robin type and Navier type boundary conditions in spaces $L^q_\si(\Om)$ and $\wt{L}^q_\si(\Om)$.  

Invariance of $L^q_\si(\Om)$ under the semigroup generated by the Hodge Laplacian on certain Helmholtz domains is used in several papers, we mention \cite{Abe:cyl}, \cite{AKST}, \cite{hobus-saal}, \cite{KUhl:ell-sys}, \cite{MM:HNS-nonlin}, \cite{MM:TAMS}.
In \cite{hobus-saal} this is even shown for some uniform $C^{2,1}$-domains without an $L^q$-Helmholtz decomposition, under the additional condition \cite[Assumption~2.4]{hobus-saal} that holds, e.g., for perturbed cones and $(\eps,\infty)$-domains (see \cite[Section~12]{hobus-saal}), but fails for aperture domains (we refer to Remark~\ref{rem:disc-As2.4} below). 
In this paper, we find $L^q$-spaces of solenoidal vector fields that are invariant for $q\in[1,\infty]$ on \emph{any} 
uniform $C^{2,1}$-domain and show invariance of the usual space $L^q_\si(\Om)$ for all $q\in(1,\infty)$ if $d=2$ and for $q\in(1,\frac{d}{d-1}\cup[2,\infty)$ in general dimension $d\ge3$.

Boundedness of the $H^\infty$-calculus in $L^q(\Om)^d$, $q\in(1,\infty)$, for the Hodge Laplacian is shown in \cite{GHT} in uniform $C^3$-domains $\Om\subseteq\R^d$. This result is used in \cite{Abe:cyl} on a cylindrical domain in $\R^3$ to show inclusion into $W^{1,q}$ of the domain of the square root. Here we show that the Hodge Laplacian enjoys a better H\"ormander functional calculus on general uniform $C^{2,1}$-domains and determine fractional domain spaces exactly (see Corollary~\ref{cor:HLap-frac-dom}). Invariance of solenoidal $L^q$-spaces then yields a H\"ormander functional calculus and, in particular, a bounded $H^\infty$-calculus for the correponding Hodge Stokes operators (see Theorem~\ref{thm:HS-Hinfty}). This in turn leads to precise descriptions the fractional domain spaces of these Hodge Stokes operators if $L^q_\si(\Om)$ is invariant under the Hodge Laplace semigroup (see Corollary~\ref{cor:HStokes-frac-dom}).

We give an overview of the paper. In Section~\ref{sec:prelim} we gathered preliminary material on boundary conditions, regularity of domains, function spaces, Helmholtz decompositions, maximal $L^p$-regularity, and functional calculi. 

In Section~\ref{sec:Hodge-Lap} we study the Hodge Laplacian on uniform $C^{2,1}$-domains. We define the operator in $L^2(\Om)^d$ by a suitable symmetric sesquilinear form and show that this coincides with the Laplacian with perfect slip boundary conditions in \cite{hobus-saal}, see Proposition~\ref{prop:H-PS}. By Davies' method we establish kernel bounds of Gaussian type for the semigroup, see Theorem~\ref{thm:Gaussian-bds}. The approach is similar to what has been done in \cite{MM:TAMS} and \cite{KUhl:ell-sys}, but we can use the precise domain descriptions from \cite{hobus-saal} to cover the full range of $q$ up to $\infty$. Then the results from \cite{DOS} or \cite{KUhl:spec-mult} apply and yield a H\"ormander functional calculus and boundedness of the $H^\infty$-calculus for the Hodge Laplacian, see Theorem~\ref{thm:HLap-Hinfty}. We also identify fractional domain spaces, see Corollary~\ref{cor:HLap-frac-dom}.

In Section~\ref{sec:Hodge-Stokes} we introduce several subspaces of solenoidal vector fields and establish invariance properties under the semigroup generated by the Hodge Laplacian, see Proposition~\ref{prop:inv-Lq}. This allows to define Hodge Stokes operators and we obtain precise domain descriptions in Proposition~\ref{prop:dom-H-Stokes}, functional calculi in Theorem~\ref{thm:HS-Hinfty}, and can identify fractional domain spaces in Corollary~\ref{cor:HStokes-frac-dom}.

In Section~\ref{sec:pert-HBC} we study Stokes operators with Robin type boundary conditions as perturbations of Hodge Stokes operators. To this end we need estimates on the solutions of the resolvent problem for the Hodge Stokes operator with inhomogeneous boundary conditions. As we dispense with \cite[Assumption 2.4]{hobus-saal} and only work under the weaker assumption that $L^q_\si(\Om)$ is invariant under the semigroup generated by the Hodge Laplacian, we reprove in Theorem~\ref{thm:Stokes-res} the resolvent estimates we need under the Assumption~\ref{ass:HH-dec}, which is familiar from \cite{hobus-saal}. Since the perturbation is of lower order, we obtain boundedness of the $H^\infty$-calculus and information on fractional domain spaces, see Theorem~\ref{thm:Robin-Stokes-Lip} and Corollary~\ref{cor:Robin-Stokes}. Similar results hold on the spaces $\wt{L}^q_\si(\Om)$ for all $q\in(1,\infty)$ without further assumptions, see Theorem~\ref{thm:Robin-Stokes-Lip-tilde} and Corollary~\ref{cor:Robin-Stokes-tilde}, but we omit the similar proofs.

We have gathered several auxiliary results in an appendix. 

Finally, we want to draw attention to the following aspect of our work. The main result of \cite{GHHS} showed that, for a uniform $C^3$-domain $\Om\subseteq\R^d$, existence of the Helmholtz projection in $L^q(\Om)^d$ (``weak Neumann'') implied maximal $L^p$-regularity, $1<p<\infty$, for the Stokes operator with dirichlet or ``no slip'' boundary conditions on $\Om$. This had been upgraded to boundedness of the $H^\infty$-calculus in \cite{GK}. Our results in this paper demonstrate in particular that ``weak Neumann'' also implies a bounded $H^\infty$-calculus for the corresponding Stokes operators with Hodge, Navier type, and Robin boundary conditions. But our results also cover certain non-Helmholtz domains.

\subsection*{Notation} As usual we understand partial derivatives $\partial_j=\frac{\partial}{\partial x_j}$, the gradient $\nabla$, the divergence operator $\dive$, or the Laplacian $\Delta$ acting on $L^1_{\mathrm{loc}}(\Om)$-functions in the distributional sense. Without explicit mentioning, we understand functions on the boundary $\partial\Om$ in the sense of traces even when we write $\ldots|_{\partial\Om}$ occasionally. We refer to the appendix for results on traces.

Sometimes, we write $a\lesssim b$ if $a\le C b$ for some inessential constant $C>0$. 

\section{Preliminaries}\label{sec:prelim}

\subsection{Boundary conditions}\label{subsec:BC}
We recall the deformation tensor 
$$
 D(u)=\frac12(\nabla u+\nabla u^T)=\frac12(\partial_ju_k+\partial_ku_j)_{j,k=1}^d
$$ 
for a vector field $u$ on subsets of $\R^d$, where we denote $\nabla u=(\partial_ju_k)_{j,k=1}^d$, i.e. $\nabla u$ has columns $\nabla u_k$. 
As in \cite{hobus-saal} we shall also use $D_\pm(u):=(\nabla u\pm\nabla u^T)$. Observe that $D_+(u)=2D(u)$, that the definition of $D_-(u)$ in \cite{hobus-saal} has the other sign, and that
$$
 \nu\times\rot u=D_-(u)\nu
$$
in case $d=3$.
We also recall the Cauchy stress tensor $T(u,p)=2D(u)-pI$, where $I\in\C^{d\times d}$ is the identity matrix and $p$ denotes the pressure. 

The boundary conditions studied in \cite{FarRost-Res} for a domain $\Om\subseteq\R^d$ with outer unit normal $\nu$ and a sufficiently smooth vector field $u$ on $\Om$ are of the form
\begin{equation}\label{eq:BC-alpha-beta}
 \nu\cdot u = 0, \qquad \alpha u +\beta[T(u,p)\nu]_{\mathrm{tan}} = 0 \quad\mbox{on $\partial\Om$},
\end{equation}
where $[\ldots]_{\mathrm{tan}}$ denotes the tangential part, and $\alpha\in[0,1)$ and $\beta\in(0,1]$ satisfy $\alpha+\beta=1$. The first condition means that the motion at the boundary is only possible in tangential directions which is reasonable for a fixed domain. Since the pressure is scalar-valued we thus have
$$
 [T(u,p)\nu]_{\mathrm{tan}}=[D_+(u)\nu]_{\mathrm{tan}},
$$
and $\alpha=0$ corresponds to Navier's slip condition where there are no tangential stress on the fluid at the boundary. The case $\beta=0$ would correspond to no-slip or Dirichlet conditions but this is excluded here. For $\alpha,\beta\in(0,1)$ one has partial slip conditions where the tangential stress at the boundary is proportional to the velocity $[u]_{\mathrm{tan}}=u$ (recall $\nu\cdot u=0$). 

In addition to these conditions, \cite{hobus-saal} also covers the conditions
\begin{equation}\label{eq:HBC}
  \nu\cdot u = 0, \qquad D_-(u)\nu = 0 \quad\mbox{on $\partial\Om$},
\end{equation}
termed ``perfect slip'' there and ``perfect wall'' in \cite{AKST}. For $d=3$ this reads
$$
  \nu\cdot u=0, \qquad \nu\times\rot u=0 \quad\mbox{on $\partial\Om$},
$$
meaning that vorticy has to be in normal direction at the boundary. This condition is called ``Hodge boundary condition'' in, e.g., \cite{MM:HNS-nonlin} since resolvents of the corresponding Laplace operator respect the Hodge (or Helmholtz) decomposition of $L^q$-vector fields on bounded Lipschitz domains, at least for an interval around $q=2$. As we shall see this partly persists also to domains $\Om$ without a Helmholtz decomposition in $L^q(\Om)^d$.
 
We refer to \cite[Section 2]{MM:HNS-nonlin} for a proof of the following fact: If the boundary of $\Om$ is of class $C^2$ and $\nu\cdot u=0$ on $\partial\Om$ then
\begin{equation}\label{eq:weingarten}
 [D_+(u)\nu]_{\mathrm{tan}}=-D_-(u)\nu+2\mathcal{W}u\quad\mbox{on $\partial\Om$},
\end{equation}
where $\mathcal{W}$ denotes the Weingarten map on $\partial\Om$ and we consider $\mathcal{W}$ as a $d\times d$-matrix-valued function on $\partial\Om$, which has values in the real symmetric matrices. In \cite{MM:HNS-nonlin} this is shown for a bounded $C^2$-domain, but the property can clearly be localized. If $\partial\Om$ is uniformly $C^2$ then $\mathcal{W}$ is continuous and bounded. If $\Om$ is a uniform $C^{2,1}$-domain then 
$\mathcal{W}$ is 
Lipschitz continuous and bounded on $\partial\Omega$. There is also a clear relation of \eqref{eq:weingarten} to \cite[Lemma 9.1]{hobus-saal}.
In view of the boundary conditions studied in \cite{hobus-saal} we remark that 
$$
 [D_-(u)\nu]_{\mathrm{tan}}= D_-(u)\nu,
$$
which via \eqref{eq:weingarten} also reflects that $\mathcal{W}u$ is tangential to $\partial\Om$ if $\nu\cdot u=0$ on $\partial\Om$.

We thus see that it is sufficient to study the Stokes operator with Hodge boundary conditions \eqref{eq:HBC} and then consider zero order perturbations of the form
\begin{equation}\label{eq:pert BC}
 \nu\cdot u=0, \qquad D_-(u)\nu = [Bu]_{\mathrm{tan}} \quad\mbox{on $\partial\Omega$},
\end{equation}
where $B\in C^{0,1}(\partial\Om)^{d\times d}$ is real-valued and symmetric. For the boundary conditions in \eqref{eq:BC-alpha-beta} we may take
\begin{align}\label{eq:Navier-Robin-BC}
 B=\al I +2 \beta \mathcal{W},
\end{align}
where $I\in\C^{d\times d}$ denotes the identity. The conditions \eqref{eq:pert BC} are called Robin boundary conditions in, e.g., \cite{Mon-Ou}. For the Weingarten map we do not need $[\ldots]_{\mathrm{tan}}$ here, but in the general case we have to put it since the left hand side in the second condition in \eqref{eq:pert BC} is tangential. 

\subsection{Regularity of domains}

We start by recalling the definition of uniform $C^{k,1}$-domains and uniform $C^k$-domains for $k\in\N_0$ and $k\in\N$, respectively.

\begin{definition}\label{def:c11}\rm
A domain $\Om\subseteq\R^d$ is called a \emph{uniform $C^{k,1}$-domain} with $k\in\N_0$  (or \emph{uniform $C^k$-domain} with $k\in\N$, respectively) 
if there are constants $\al,\beta,K>0$ such that, for each $x_0\in\partial\Om$,
there is a Cartesian coordinate system with origin at $x_0$ and coordinates $y=(y',y_d)$,
$y'=(y_1,\ldots,y_{d-1})$ and a $C^{k,1}$-function (or $C^k$-function, respectively) $h$, defined on $\{y':|y'|\le\al\}$
and with $\|h\|_{C^{k,1}}\le K$ (or $\|h\|_{C^k}$, respectively), such that, for the neighborhood 
$$
 U_{\al,\beta,h}(x_0)=\{y=(y',y_d)\in\R^d:|y_d-h(y')|<\beta,|y'|<\al\}
$$
of $x_0$, we have $U_{\al,\beta,h}(x_0)\cap\partial\Om=\{(y',h(y')):|y'|<\al\}$ and
$$
 U_{\al,\beta,h}(x_0)\cap\Om=\{(y',y_d): h(y')-\beta<y_n<h(y'),|y'|<\al\}.
$$ 
A domain $\Om\subseteq \R^d$ is called a \emph{Lipschitz domain} if, for every $x_0\in\partial\Om$, one can find $\al,\beta, K>0$ and a Lipschitz function $h$ such that one has a representation as above, and $\Om$ is called a \emph{uniform Lipschitz domain} if $\Omega$ is a uniform $C^{0,1}$-domain. 
Observe that for a Lipschitz domain, the constants $\alpha,\beta,K>0$ in the local representation may depend on $x_0\in\partial\Om$.   
\end{definition}

\subsection{Function spaces and Helmholtz decompositions}\label{subsec:function-spaces}
Let $\Om\subseteq\R^d$ be a domain. For $q\in[1,\infty]$ and $k\in\N$, $L^q(\Om)$ and $W^{k,q}(\Om)$ denote the usual Lebesgue and Sobolev spaces on $\Om$: 
$$
 W^{k,q}(\Om):=\{u\in L^q(\Om):\partial^\alpha u\in L^q(\Om) \ \mbox{for all $\alpha\in\N_0^d$ with $|\alpha|\le k$}\,\}.
$$
Without explicit further notice functions in these spaces are always complex-valued.

We denote by $C_c(\ov{\Om})$ the space of continuous functions on $\ov{\Om}$ that have compact support in $\ov{\Om}$ and
by $C_0(\ov{\Om})$ the closure of $C_c(\ov{\Om})$ w.r.t. the sup-norm $\nn{\cdot}{\infty}$ in, e.g., the space $C_b(\ov{\Om})$ of all bounded and continuous functions on $\ov{\Om}$ or in $L^\infty(\Om)$.

We denote by $C^\infty_c(\Om)$ the set of all $C^\infty$-functions on $\Om$ with compact support in $\Om$. Then we denote $C^\infty_{c,\si}(\Om):=\{\ph\in C^\infty_c(\Om)^d: \dv\ph=0\}$ and, for $q\in(1,\infty)$, 
$L^q_\si(\Om)$ denotes the closure of $C^\infty_{c,\si}(\Om)$ in $L^q(\Om)^d$.

We denote by $G^q(\Om)$ the space of $L^q$-gradients, i.e., the space of all $f\in L^q(\Om)^d$ such there exists a distribution $\psi$ on $\Om$ with $f=\nabla\psi$. It is well-known that we then have $\psi\in L^q_\loc(\Om)$, i.e. $\psi|_K\in L^q(K)$ for any compact subset $K\subseteq\Om$, and that, for a Lipschitz domain $\Om\subseteq\R^d$, we even have $\psi\in L^q_\loc(\ov{\Om})$, i.e. $\psi|_K\in L^q(K)$ for any compact subset $K\subseteq\ov{\Om}$.

We thus set, for a Lipschitz domain $\Om\subseteq\R^d$,
$$
 \wh{W}^{1,q}(\Om):=\{\psi\in L^q_\loc(\ov{\Om}) : \nabla \psi\in L^q(\Om)^d\}/\C,
$$ 
so that $G^q(\Om)=\nabla\wh{W}^{1,q}(\Om)$.

For the usual duality between $L^q(\Om)^d$ and $L^{q'}(\Om)^d$, where $q'\in(1,\infty)$ denotes the dual exponent to $q$ given by $\frac1q+\frac1{q'}=1$, we then have
\begin{align}\label{eq:perp-spaces}
 L^q_\si(\Om)^\perp=G^{q'}(\Om), \qquad G^q(\Om)^\perp=L^{q'}_\si(\Om).
\end{align}
For $q=2$ one has the orthogonal decomposition $L^2(\Om)^d=L^2_\si(\Om)\oplus G^2(\Om)$, usually called \emph{Helmholtz} or (a special case of) \emph{Hodge decomposition}, with corresponding orthogonal projection $\PP_2$ in $L^2(\Om)^d$, the  \emph{Helmholtz projection}.

\begin{remark}\label{rem:cons-HHdec}\rm
Let $\Om\subseteq\R^d$ be a uniform $C^{1}$-domain. It is shown in \cite[Proposition 2.1]{GK} that
there is an interval $I_\PP\subseteq(1,\infty)$ with $2\in I_\PP$ and symmetric in the sense that $q\in I_\PP$ if and only if $q'\in I_\PP$, such that, for any $q\in I_\PP$, $\PP_2$ restricted to $L^2(\Om)^d\cap L^q(\Om)^d$ extends to a bounded operator $\PP_q$ on $L^q(\Om)^d$ related to the Helmholtz decomposition $L^q(\Om)^d=L^q_\si(\Om)\oplus G^q(\Om)$. For $q\in I_\PP$, $\PP_q$ is called \emph{Helmholtz projection} in $L^q(\Om)^d$.
Moreover, for $q\in(1,\infty)$ we have $q\in I_\PP$ if and only if the Helmholtz decomposition $L^q(\Om)^d=L^q_\si(\Om)\oplus G^q(\Om)$ holds.
\end{remark}

It is well-known that, for a bounded Lipschitz domain $\Om\subseteq\R^d$, one can give a sense to the normal component $\nu\cdot f$ on the boundary $\partial\Om$ for vector-fields $f\in L^q(\Om)^d$ with $\dv f\in L^q(\Om)$ (we refer, e.g., to \cite[II.1.2]{Sohr}). This is done via an integration by parts formula and, since it can be localized, persists to uniform Lipschitz domains (see also Proposition~\ref{prop:app-trace} in the Appendix for details in uniform $C^{2,1}$-domains). 
For a bounded Lipschitz domain one has
$$
 L^q_\si(\Om)=\{f\in L^q(\Om)^d: \dv f=0, \nu\cdot f|_{\partial\Om}=0\,\}.
$$
For a general Lipschitz domain $\Om\subseteq\R^d$ we denote 
$$
 \mathcal{L}^q_\si(\Om)=\{f\in L^q(\Om)^d: \dv u=0, \nu\cdot u|_{\partial\Om}=0\,\},\qquad
 \mathcal{G}^q(\Om)=\ov{\nabla C^\infty_c(\ov{\Om})}^{L^q(\Om)^d}.
$$
Clearly, $L^q_\si(\Om)\subseteq\mathcal{L}^q_\si(\Om)$ and $\mathcal{G}^q(\Om)\subseteq G^q(\Om)$, and there are unbounded domains with several outlets to infinity where those inclusions are strict with arbitrary finite or with infinite codimension, see \cite{Bog86}, \cite{MaBo}. We have the duality relations
\begin{align}\label{eq:perp-calspaces}
 \mathcal{L}^q_\si(\Om)^\perp=\mathcal{G}^{q'}(\Om), \qquad \mathcal{G}^q(\Om)^\perp=\mathcal{L}^{q'}_\si(\Om),
\end{align}
and for $q=2$ the orthogonal decomposition $L^2(\Om)^d=\mathcal{L}^2_\si(\Om)\oplus \mathcal{G}^2(\Om)$ with corresponding orthogonal projection $\mathcal{P}_2$ in $L^2(\Om)^d$.

\begin{remark}\rm\label{rem:disc-As2.4}
We discuss \cite[Assumption 2.4]{hobus-saal}, essential for a number of results in \cite{hobus-saal}, some of which we shall extend. This assumption reads: $\nabla C^\infty_c(\ov{\Om})$ is dense in $G^{q'}(\Om)$. By \eqref{eq:perp-spaces} this is equivalent to $L^q_\si(\Om)=\big(\nabla C^\infty_c(\ov{\Om})\big)^\perp=\mathcal{G}^{q'}(\Omega)^\perp$. 
By \eqref{eq:perp-calspaces} it is finally equivalent to $L^q_\si(\Om)=\mathcal{L}^q(\Omega)$. 

This always holds for $q\in(1,\frac{d}{d-1}]$, see Lemma~\ref{lem:Lq-solenoidal} below, but for large $q$ it fails, e.g., for aperture domains or other domains with several outlets to infinity, see \cite{mas-bog-appr}.
\end{remark}

As for the Helmholtz decomposition above, there is an interval $I_\mathcal{P}\ni 2$, symmetric in the sense that $q\in I_\mathcal{P}$ if and only if $q'\in I_\mathcal{P}$, such that, for any $q\in I_\mathcal{P}$, $\mathcal{P}_2$ restricted to $L^2(\Om)^d\cap L^q(\Om)^d$ extends to a bounded projection $\mathcal{P}_q$ on $L^q(\Om)^d$ related to the decomposition $L^q(\Om)^d=\mathcal{L}^q_\si(\Om)\oplus \mathcal{G}^q(\Om)$. 
The proofis very similar to the prooof of \cite[Proposition 2.1]{GK}.

As we shall also use results from \cite{FKS:ArchMath} on the Helmholtz condition on uniform $C^1$-domains we recall the spaces
$$
 \wt{L}^q(\Om):=\left\{\begin{array}{ll}
                           L^q(\Om)\cap L^2(\Om), & 2\le q<\infty,\\
                            L^q(\Om)+ L^2(\Om), & 1<q<2,\\
                           \end{array}
                           \right.
$$
and the corresponding spaces of solenoidal and gradient vector fields
$$
 \wt{L}^q_\si(\Om):=\left\{\begin{array}{ll}
                           L^q_\si(\Om)\cap L_\si^2(\Om), & 2\le q<\infty,\\
                            L^q_\si(\Om)+ L_\si^2(\Om), & 1<q<2,\\
                           \end{array}
                           \right.
\qquad
 \wt{G}^q(\Om):=\left\{\begin{array}{ll}
                           G^q(\Om)\cap G^2(\Om), & 2\le q<\infty,\\
                            G^q(\Om)+ G^2(\Om), & 1<q<2.\\
                           \end{array}
                           \right.
$$
For $k\in\N$ and $q\in(1,\infty)$, we shall later on also meet the spaces
$$
 \wt{W}^{k,q}(\Om):=\left\{\begin{array}{ll}
                           W^{k,q}(\Om)\cap W^{k,2}(\Om), & 2\le q<\infty,\\
                            W^{k,q}(\Om)+ W^{k,2}(\Om), & 1<q<2.\\
                           \end{array}
                           \right.
$$
We state explicitly that, in the usual canonical way,
\begin{equation}\label{eq:tilde-Lq-dual}
 \big(\wt{L}^q(\Om)\big)'=\wt{L}^{q'}(\Om),\qquad 1<q<\infty.
\end{equation}
The following on the Helmholtz decomposition in $\wt{L}^q(\Om)^d$ has been shown in \cite[Theorem 1.2, Corollary 1.3]{FKS:ArchMath}. We remark that this has been used in the proof of the assertion of Remark~\ref{rem:cons-HHdec} in \cite[Proposition~2.1]{GK}.

\begin{theorem}\label{thm:FKS-dec}
Let $\Om\subseteq\R^d$ be a uniform $C^1$-domain and $1<q<\infty$. Then
$$
 \wt{L}^q(\Om)^d=\wt{L}_\si(\Om)\oplus\wt{G}^q(\Om)
$$
and the correponding projection $\wt{P}_q$ in $\wt{L}^q(\Om)^d$ satisfies $(\wt{P}_q)'=\wt{P}_{q'}$.

Moreover, $C^\infty_{c,\si}(\Om)$ is dense in $\wt{L}^q_\si(\Om)$ for the norm of $\wt{L}^q(\Om)^d$ and one has the annihilator relations
$$
  \wt{L}^q_\si(\Om)^\perp=\wt{G}^{q'}(\Om), \qquad \wt{G}^q(\Om)^\perp=\wt{L}^{q'}_\si(\Om),
$$
and in a canonical way the isomorphisms
$$
 \big(\wt{L}^q_\si(\Om)\big)'\simeq \wt{L}^{q'}_\si(\Om), \qquad
 \big(\wt{G}^q(\Om)\big)' \simeq \wt{G}^{q'}(\Om).
$$ 
\end{theorem}

\subsection{Maximal $L^p$-regularity, $H^\infty$-functional calculus, and H\"ormander functional calculus}\label{subsec:Hinfty} 
We only recall basic notions and refer to \cite{KuW:levico} for more details. Let $-A$ be the densely defined generator of a bounded analytic semigroup in a Banach space $X$. For $p\in(1,\infty)$, $A$ is said to have \emph{maximal $L^p$-regularity} if, for any $f\in L^p(\R_+;X)$, there exists a unique mild solution of the Cauchy problem
$$
 u'(t)+Au(t)=f(t),\quad t>0,\qquad u(0)=0,
$$
which satisfies $u',Au\in L^p(\R_+;X)$. The densely defined negative generator $B$ of an analytic semigroup is said to have 
\emph{maximal $L^p$-regularity on finite intervals} if, for some (and then equivalently for all) $T>0$ and any $f\in L^p(0,T;X)$, there exists a unique mild solution of the Cauchy problem
$$
 u'(t)+Bu(t)=f(t),\quad t\in(0,T),\qquad u(0)=0,
$$
which satisfies $u',Bu\in L^p(0,T;X)$. If $A$ has maximal $L^p$-regularity then any translate $B=\mu+A$, $\mu\in\R$, has maximal $L^p$-regularity on finite intervals. Conversely, if $B$ has maximal $L^p$-regularity on finite intervals then $B+\mu$ has maximal $L^p$-regularity for some $\mu\ge0$.  

In UMD spaces $X$, in particular in closed subspaces of $L^q$-spaces with $q\in(1,\infty)$, maximal $L^p$-regularity for $p\in(1,\infty)$ is characterized by $R$-sectoriality of $A$ of some angle $<\frac\pi2$ (see, e.g., \cite[1.11]{KuW:levico}). Here, the operator $A$ is called \emph{$R$-sectorial of angle $\om\in[0,\pi)$} if $\si(A)\subseteq\Si_\om:=\{\la\in\C\sm\{0\}:|\arg\la|\le\om\}\cup\{0\}$, and for any $\theta\in(\om,\pi)$, the set $\{\la R(\la,A):\la\in\C\sm\Si_\theta\}\subseteq\cL(X)$ is $R$-bounded.

For Banach spaces $X$, $Y$ a subset $\tau\subseteq\cL(X,Y)$ is called \emph{$R$-bounded} with \emph{$R$-bound} $C$ if, for all $n\in\N$, $x_1,\ldots,x_n\in X$ and $T_1,\ldots,T_n\in\tau$ one has
$$
 \EE\big\|\sum_{j=1}^n\eps_jT_jx_j\big\|_{Y}\le C\EE\big\|\sum_{j=1}^n\eps_jx_j\big\|_{X}, 
$$
where the $\eps_j$ are independent and symmetric $\{-1,1\}$-valued random variables, e.g., Rachemachers. By the Khintchine-Kahane inequalities, for $X=L^q$ with $q\in(1,\infty)$, expressions $\EE\n{\sum_j\eps_jf_j}$ are equivalent to square function expressions $\n{(\sum_j|f_j|^2)^{1/2}}$. This has been extensively used in, e.g., \cite{BK-mr}.

If we replace, in the above definition of $R$-sectoriality, $R$-boundedness by boundedness, we obtain the definition of a \emph{sectorial operator of angle $\om\in[0,\pi)$}. A sectorial operator $A$ of angle $\om\in[0,\pi)$ is said to have a bounded
$H^\infty(\Si_\theta)$-calculus, where $\theta\in(\om,\pi)$, if for some $C>0$ we have the bound
$$
 \n{F(A)}\le C\nn{F}{\infty,\Si_\theta}
$$
for all $F$ holomorphic on the interior of $\Si_\theta$, for which $|F(z)|\le M\min\{|z|^\eps,|z|^{-\eps}\}$ holds for some $M,\eps>0$. Here, the operator $F(A)\in\cL(X)$ is defined by the Cauchy type integral
\begin{equation}\label{eq:Hinfty-contour}
 F(A)=\frac1{2\pi i}\int_{\partial\Si_\si} F(\la) R(\la,A)\,d\la,
\end{equation}
with $\si\in(\om,\theta)$. Observe that this is a Bochner integral by the assumptions on $F$.

If $A$ is densely defined with dense range and has a bounded $H^\infty(\Si_\theta)$-calculus then $F(A)$ is a bounded operator for all $F$ holomorphic and bounded on the interior of $\Si_\theta$. In particular, $A$ has fractional powers $A^{it}\in\cL(X)$ for all $t\in\R$, with an exponential bound in $|t|$, i.e. \emph{$A$ has bounded imaginary powers}. It is well-known that, if $A$ has bounded imaginary powers, then for $\theta\in(0,1)$ the domains fo the fractional powers $A^\theta$ are obtained by complex interpolation
$$
 D(A^\theta)=[X,D(A)]_\theta,\qquad\theta\in(0,1),
$$  
see, e.g., \cite{Triebel}, \cite{KuW:levico}, \cite{fc-book}.

Under the same assumptions, the operator $A$ is said to have a \emph{H\"ormander functional calculus} if there exist $C>0$ and $s>0$ such that, for some $\eta\in C_c^\infty(0,\infty)\sm\{0\}$, one has an estimate
\begin{equation}\label{eq:H-Ws2-cond}
 \n{F(A)}\le C\sup_{t>0}\nn{\eta(\cdot)F(t\cdot)}{W^{s,2}}
\end{equation}
for $F\in C^\infty_c(0,\infty)$, say. For more on this type of functional calculus we refer to \cite{DOS}, \cite{KUhl:spec-mult}, \cite{KUhl:ell-sys}. In the typical situation $X=L^q(\Om)$, $A$ is self-adjoint in $L^2(\Om)$ and, at least on $L^q(\Om)\cap L^2(\Om)$, the operator $F(A)$ is given by the spectral theorem in $L^2(\Om)$. Let us already mention here that we do not aim for optimality of the smoothness parameter $s$ here and view this property more as a qualitative strengthening of a bounded $H^\infty$-calculus: If $F$ is bounded and holomorphic on the interior of $\Si_\theta$ then, for any $t>0$ and $k\in\N$, we have by Cauchy's integral formula 
\begin{align*}
 F^{(k)}(t)=\frac{k!}{2\pi i}\int_{|z-t|\le ct}\frac{F(z)}{(z-t)^{k+1}}\,dz, 
\end{align*}
for any $c\in(0,\arcsin\theta)$, which leads to $|t^k F^{(k)}(t)|\le \frac{k!}{c^k}\nn{F}{\infty,\Si_\theta}$.
This shows that a H\"ormander functional calculus for some $s>0$ implies a bounded $H^\infty(\Si_\theta)$-calculus for any $\theta\in(0,\frac\pi2)$.

\section{The Hodge Laplacian on uniform $C^{2,1}$-domains}\label{sec:Hodge-Lap}

In this section we study the so-called Hodge Laplacian in uniform $C^{2,1}$-domains. We establish pointwise Gaussian kernel bounds for the semigroup operators. Similar to the approach in \cite{KUhl:ell-sys} this is done by Davies' method. Compared to the situation in bounded Lipschitz domains in \cite{KUhl:ell-sys} we can here make use of the $L^q$-theory of \cite{hobus-saal}, in particular the precise description of the domain of the operator in $L^q(\Om)^d$, and combine this with Sobolev embeddings. An application of the main result of \cite{DOS} then yields a bounded H\"ormander functional calculus on the $L^q$-scale. This calculus is much stronger than a bounded $H^\infty$-calculus, for which an application of the main result in \cite{duong-robinson} would have been sufficient. In any case this leads to bounded imaginary powers and thus to a precise description of the domain of the square root of the operator in $L^q(\Om)^d$.     

\subsection{The operator}
We define the Hodge Laplacian $\Delta_H$ in $L^2(\Om)^d$ for a Lipschitz domain $\Om\subseteq\R^d$ by a suitable sesquilinear form. For $d=3$ we recall the following from \cite{MM:HNS-nonlin}, \cite{KUhl:ell-sys}. Let
\begin{equation}\label{eq:form-rot}
 \a :V\times V\to \C,\qquad \a (u,v):=\int_\Om \rot u\cdot\ov{\rot v}\,dx + \int_\Om \dv u\,\ov{\dv v}\,dx, 
\end{equation}
where 
$$
   V:=V(\Om):=\{u\in L^2(\Om)^3: \rot u\in L^2(\Om)^3,\ \dv u\in L^2(\Om),\ \nu\cdot u|_{\partial\Om}=0\}.
$$
Notice that the boundary condition in the definition of $V$ makes sense.
Then $-\Delta_H$ is the operator in $L^2(\Om)^3$ associated with $\a $ in the usual sense: For $u,f\in L^2(\Om)^3$
we have $u\in D(\Delta_H)$ and $-\Delta_H u=f$ if and only if 
$$
 u\in V\qquad\mbox{and}\qquad \forall v\in V: \a (u,v)=\sk{f}{v},
$$
where $\sk{f}{v}=\int_\Om f\cdot\ov{v}\,dx$ denotes the scalar product in $L^2(\Om)^3$. 
For $d\ge3$ we take inspiration from \cite{MM:TAMS} (see also the weak formulation in \cite{AKST}) and let
\begin{equation}\label{eq:form-a}
 \a :V\times V\to \C, \qquad  
 \a (u,v):=\frac12\int_\Om D_-(u):\ov{D_-(v)}\,dx + \int_\Om \dv u\,\ov{\dv v}\,dx, 
\end{equation}
where
$$
  V:=V(\Om):=\{u\in L^2(\Om)^d: D_-(u)\in L^2(\Om)^{d\times d},\ \dv u\in L^2(\Om),\ \nu\cdot u|_{\partial\Om}=0\}
$$ 
and 
$$
 B_1:\ov{B_2}:=\sum_{j,k=1}^d b^1_{jk}\ov{b^2_{jk}}\quad\mbox{for matrices $B_{l}=(b^{l}_{jk})_{jk}\in \C^{d\times d}$, $l=1,2$.}
$$

\begin{prop}
Let $\Omega\subseteq\R^d$ be a Lipschitz domain.
Then the operator $-\Delta_H$ associated with $\a $ in $L^2(\Om)^d$ is self-adjoint in $L^2(\Om)^d$ and 
$-\Delta_H\ge0$.
\end{prop}

\begin{proof}
The sesquilinear form $\a $ is symmetric, i.e. $\a (u,v)=\ov{\a (v,u)}$ for all $u,v\in V(\Om)$.
The space $V(\Om)$ is dense in $L^2(\Om)^d$ (since it contains $H^1_0(\Om)^d$) and $V$ is a Hilbert space for the scalar product
$$
 \sk{u}{v}_V:=\a (u,v)+\sk{u}{v}_{L^2(\Om)^d}.
$$
Hence the operator $-\Delta_H$ associated with $\a $ in $L^2(\Om)^d$ is self-adjoint in $L^2(\Om)^d$.
By $\a (u,u)\ge0$ for all $u\in V(\Om)$, $-\Delta_H$ is non-negative.
\end{proof}

\begin{corollary}
Let $\Om\subseteq\R^d$ be a Lipschitz domain. Then $\Delta_H$ generates a bounded analytic semigroup 
$(T(t))_{t\ge0}:=(e^{t\Delta_H})_{t\ge0}$ in $L^2(\Om)^d$ which is contractive on $\{z\in\C:\mathrm{Re}z>0\}$.
\end{corollary}

We determine the operator $-\Delta_H$ associated with $\a $, assuming additional regularity of the boundary.
To this end we also need the following result which is part of \cite[Theorem 6.1]{hobus-saal}.

\begin{prop}\label{prop:Lap-PS}
Let $\Om\subseteq\R^d$ be a uniform $C^{2,1}$-domain and $q\in(1,\infty)$. The restriction 
$\Delta_{PS,q}$ of the Laplacian $\Delta$ to the set 
$$
 D(\Delta_{PS,q}):=\{u\in W^{2,q}(\Om)^d: \nu\cdot u=0\ \mbox{and}\ D_-(u)\nu=0\ \mbox{on $\partial\Om$}\,\}
$$ 
is the generator of an analytic semigroup in $L^q(\Om)^d$.
\end{prop}

\begin{prop}\label{prop:H-PS}
Let $\Om\subseteq\R^d$ be a uniform $C^{2,1}$-domain. Then  $-\Delta_H$ coincides with the operator $-\Delta_{PS,2}$, i.e.
$$
 D(-\Delta_H)=\{u\in W^{2,2}(\Om)^d:\nu\cdot u=0\ \mbox{and}\ D_-(u)\nu=0\ \mbox{on $\partial\Om$}\,\}
$$
and, for $u\in D(-\Delta_H)$,
$$
 -\Delta_H u=-\Delta u.
$$
Moreover we have
$$
 V(\Om)=\{u\in W^{1,2}(\Om)^d: \nu\cdot u=0\ \mbox{on $\partial\Om$}\,\}.
$$
\end{prop}

\begin{corollary}\label{cor:PS-sa}
Let $\Om\subseteq\R^d$ be a uniform $C^{2,1}$-domain. Then $\Delta_{PS,2}$ is self-adjoint in $L^2(\Om)^d$ and $-\Delta_{PS,2}\ge0$.
\end{corollary}

\begin{proof}[Proof of Proposition~\ref{prop:H-PS}]
We start with the elementary formula (see also \cite[Lemma 5.3 (i)]{hobus-saal})
$$
 \dv(D_-(u)\ov{v})=(\Delta u-\nabla\dv u)\cdot\ov{v}+D_-(u):\ov{\nabla v}
$$
and symmetrize the second term
\begin{align}\label{eq:symm-D-}
 D_-(u):\ov{\nabla v}=D_-(u)^T:\ov{\nabla v}^T=-D_-(u):\ov{\nabla v}^T=\frac12 D_-(u):\ov{D_-(v)}
\end{align}
to arrive at
\begin{align}\label{eq:formula-D-}
 \frac12 D_-(u):\ov{D_-(v)}= \dv(D_-(u)\ov{v})-(\Delta u-\nabla\dv u)\cdot\ov{v}.
\end{align}
Then we use Gau{\ss}' theorem (see Proposition~\ref{prop:app-trace}) and obtain, for $u\in W^{2,2}(\Om)^d\cap V(\Om)$ and $v\in V(\Om)$,
\beq
 \a (u,v)&=&\int_\Om \dv(D_-(u)\ov{v})\,dx-\int_\Om(\Delta u-\nabla\dv u)\cdot\ov{v}\,dx
                               +\int_\Om \dv((\dv u)\ov{v})-(\nabla\dv u)\cdot\ov{v}\,dx\\
 &=&\int_\Om (-\Delta u)\cdot\ov{v}\,dx+\int_\Om \dv(D_-(u)\ov{v})\,dx+\int_\Om \dv((\dv u)\ov{v})\,dx\\
 &=&\int_\Om (-\Delta u)\cdot\ov{v}\,dx+\int_{\partial\Om}\nu\cdot D_-(u)\ov{v}\,d\si
     +\int_{\partial\Om}(\dv u)(\nu\cdot\ov{v})\,d\si\\
 &=&\int_\Om (-\Delta u)\cdot\ov{v}\,dx-\int_{\partial\Om}\ov{v}\cdot D_-(u)\nu\,d\si.
\eeq
In the last step we used $\nu\cdot D_-(u)\ov{v}=-\ov{v}\cdot D_-(u)\nu$ (see also \cite[Lemma 5.3 (iii)]{hobus-saal}) 
and $\nu\cdot v=0$ on $\partial\Om$ by $v\in V(\Om)$. This shows $-\Delta_Hu=-\Delta u$ if $u\in W^{2,2}(\Om)\cap V(\Om)$ satisfies in addition $D_-(u)\nu=0$.

Observing $W^{2,2}(\Om)^d\cap V(\Om)=\{u\in W^{2,2}(\Om): \nu\cdot u|_{\partial\Om}=0\}$ we thus have shown
$$
 D(\Delta_{PS,2})=\{u\in W^{2,2}(\Om)^d: \nu\cdot u=0\ \mbox{and}\ D_-(u)\nu=0\ \mbox{on $\partial\Om$}\,\}\subseteq D(-\Delta_H)
$$
and that $\Delta_{PS,2}$ is a restriction of $\Delta_{H}$. Since the resolvent sets of both operators $\Delta_H$ and $\Delta_{PS,2}$ contain a right half plane (here we use Proposition \ref{prop:Lap-PS}) we conclude $\Delta_H=\Delta_{PS,2}$ as claimed.

The last assertion is obtained by complex interpolation. It suffices to show $V(\Om)\subseteq W^{1,2}(\Om)$. Since $-\Delta_H$ is self-adjoint we have 
$$
 V(\Om)=[L^2(\Om)^d,D(-\Delta_H)]_{1/2}\subseteq[L^2(\Om)^d,W^{2,2}(\Om)^d]_{1/2}=W^{1,2}(\Om)^d,
$$
where we refer to Proposition~\ref{prop:app-interpol} for the last equality.
\end{proof}

For later purposes we note the following variants of the integration by parts argument in the previous proof under relaxed conditions.

\begin{lemma}\label{lem:calc-proof3.4} 
Let $\Om\subseteq\R^d$ be a uniform $C^{2,1}$-domain and $q\in(1,\infty)$.
\begin{enumerate}[label=$(\roman*)$]
\item\label{item:calc-pf3.4-i} If $u\in W^{2,q}(\Om)^d\cap L^q_\si(\Om)$ with $D_-(u)\nu=0$ on $\partial\Om$ and $v\in G^{q'}(\Om)$ then 
$$
 \int_\Om (-\Delta u)\cdot\ov{v}\,dx=0.
$$
\item\label{item:calc-pf3.4-ii} If $\nabla\psi\in G^q(\Om)$, $\Delta\psi\in W^{1,q}(\Om)$ with $\nu\cdot\nabla\psi=0$ on $\partial\Om$ and $v\in D(\Delta_{H,q'})$ then $\Delta\nabla\psi\in L^q(\Om)^d$ and 
$$
 \int_\Om\nabla\psi\cdot\ov{(-\Delta v)}\,dx=\int_\Om (-\Delta\nabla\psi)\cdot\ov{v}\,dx.
$$
\end{enumerate}
\end{lemma}

\begin{proof}
\ref{item:calc-pf3.4-i}: First observe that $D_-(v)=0$ since $v$ is a gradient. If, in addition, $v\in W^{1,q'}(\Om)^d$ then the calculations in the proof of Proposition~\ref{prop:H-PS} show the assertion. By \cite[Theorem 5]{mas-bog-appr} we can approximate a given gradient $v\in G^{q'}(\Om)$ in $L^{q'}(\Om)^d$-norm by gradients in $W^{1,q'}(\Om)$.

\ref{item:calc-pf3.4-ii} Again, we observe $D_-(\nabla\psi)=0$ in $\Om$. 
We also observe that $\Delta\nabla\psi=\nabla\Delta\psi\in L^{q}(\Om)^d$ by $\Delta\psi\in W^{1,q}(\Om)$. The formula \eqref{eq:formula-D-} still is true if just $u\in L^q(\Om)^d$ with $\dive u\in L^q(\Om)$, $D_-(u)\in W^{1,q}(\Om)^{d\times d}$, $\Delta u,\nabla \dive u\in L^q(\Om)^d$
and $v\in L^{q'}(\Om)^d$ with $\dive v\in L^{q'}(\Om)$, $D_-(v)\in L^{q'}(\Om)^{d\times d}$. 
The argument in the proof of \cite[Theorem 3.22]{adams} shows that we can approximate such a given $v$ by a sequence of smooth $v_n$ with compact support such that $v_n\to v$, $\dive v_n\to\dive v$, and $D_-(v_n)\to D_-(v)$ in $L^{q'}$-norm. Hence we can carry our the symmetrization \eqref{eq:symm-D-} for $v_n$ and pass to the limit.

Consequently we have, for $v\in D(\Delta_{H,q'})$ and $u\in L^q(\Om)^d$ with $\dive u\in L^q(\Om)$, $D_-(u)\in W^{1,q}(\Om)^{d\times d}$, and $\Delta u,\nabla \dive u\in L^q(\Om)^d$, that 
$$
 \int_\Om u\cdot\ov{(-\Delta v)}\,dx+\int_{\partial\Om}\ov{\dive v}\, (\nu\cdot u)-u\cdot \ov{D_-(v)\nu}\,d\si
=\int_\Om(-\Delta u)\cdot\ov{v}\,dx. 
$$ 
Putting $u=\nabla\psi$ as in the assumption this proves the claim. 
\end{proof}

\subsection{Gaussian bounds}
We employ the method from \cite{KUhl:ell-sys} to establish Gaussian type bounds for $(T(t))_{t\ge0}$, but here we are in a more regular situation, and can fully exploit the information in Proposition~\ref{prop:Lap-PS} on the domain of the generator in $L^q(\Om)^d$ for $2\le q<\infty$.

\begin{theorem}\label{thm:Gaussian-bds}
Let $\Om\subseteq\R^d$ be a uniform $C^{2,1}$-domain. Then the semigroup $(T(t))_{t\ge0}$ generated by $\Delta_H$ in $L^2(\Om)^d$ consists for $t>0$ of integral operator with $\R^{d\times d}$-valued integral kernels $k(t,x,y)$ satisfying pointwise Gaussian bounds, i.e., there exist constants $C,\del,b>0$ such that, for all $t>0$ and $x,y\in\Om$,
\begin{equation}\label{eq:kernel-bds}
 |k(t,x,y)|\le C t^{-d/2}\,e^{\del t}\, e^{-b\frac{|x-y|^2}t}.
\end{equation}
\end{theorem}

\begin{proof}
First we show that each $T(t)$ leaves $L^2(\Om;\R^d)$ invariant. We use \cite[Theorem 2.1]{Ouh:inv-convex}. So let $u\in V(\Om)$. We have to show $\Rep u\in V(\Om)$ which is clear and $\Rep\a(u,u-\Rep u)\ge0$.
But 
\beq
 \Rep \a(u, u-\Rep u)&=&\Rep\left(-i\int_\Om D_-(u):D_-(\Imp u)+\dv u \,\dv(\Imp u)\,dx\right)\\
 &=&\int_\Om D_-(\Imp u):D_-(\Imp u)+|\dv(\Imp u)|^2\,dx\ge0. 
\eeq
We only sketch the part of the proof in $L^2(\Om)^d$ (steps 1 and 2 below) where calculations are just as in the proof of \cite[Theorem 5.1]{KUhl:ell-sys} (there, $\Om\subseteq\R^3$ was a bounded Lipschitz domain and $\a$ given by \eqref{eq:form-rot}).
We shall use Davies' method and consider ``twisted'' forms
$$
\a _{\vrh\phi}(u,v):=\a (e^{\vrh\phi}u,e^{-\vrh\phi}v)
\qquad(u,v\in V(\Om)),
$$
where $\vrh\in\R$ and $\phi\in\E:=\{\phi\in C_c^\infty(\ov\Om,\R)
\,:\, \|\partial_j\phi\|_\infty\leq1 \mbox{ for all
}j \}$. Observe that $e^{\vrh\phi}u\in V(\Om)$ for $u\in V(\Om)$ so that $\a_{\vrh\phi}$ is well-defined.

\emph{Step 1:} For each $\ga\in(0,1)$ there exists a
constant $\om_0\geq0$ such that, for all $u\in V(\Om)$, $\vrh\in\R$,
and $\phi\in\E$,
\begin{equation}\label{eq:step1}
\bigl|\a_{\vrh\phi}(u,u)-\a(u,u)\bigr|\leq\gamma\a(u,u)+\om_0\vrh^2\|u\|_2^2\,.
\end{equation}
\emph{Step 2:} There are constants $C,\om_1>0$ such that
\begin{eqnarray}
\label{eq:2-2-est}
 \bigl\|e^{-\vrh\phi}e^{t\Delta_H}(e^{\vrh\phi}f)\bigr\|_{L^2(\Om)^d}
 &\leq& C e^{\om_1\vrh^2t}\nn{f}{L^2(\Om)^d} \\
\label{eq:D-2-est}
 \bigl\|D_-(e^{-\vrh\phi}e^{t\Delta_H}e^{\vrh\phi}f)\bigr\|_{L^2(\Om)^d}
 &\leq& C t^{-1/2} e^{\om_1\vrh^2t}\nn{f}{L^2(\Om)^d} \\
\label{eq:div-2-est}
 \bigl\|\dv (e^{-\vrh\phi} e^{t\Delta_H}e^{\vrh\phi}f)\bigr\|_{L^2(\Om)^d}
 &\leq& C t^{-1/2} e^{\om_1\vrh^2t}\nn{f}{L^2(\Om)^d}
\end{eqnarray}
for all $\vrh\in\R$, $\ph\in\E$, $t>0$ and $f\in L^2(\Om)^d$. Here we remark that, for a scalar function $g$ and a vector field $u$, we have
\begin{equation}\label{eq:mult-D-}
 D_-(g u)=g D_-(u)+(\nabla g)u^T-u(\nabla g)^T,
\end{equation}
and this formula replaces the formula $\rot(g u)=g\,\rot u+\nabla g\times u$, used in \cite{KUhl:ell-sys}.

\emph{Step 3:} We make use of the Sobolev embedding $V(\Om)\emb W^{1,2}(\Om)^d\emb L^{q_0}(\Om)^d$ 
where, for $d\ge3$, $q_0$ is given by $\frac1{q_0}=\frac12-\frac1d$, i.e. $q_0=\frac{2d}{d-2}$ (for $d=2$ see 
Remark~\ref{rem:kernel-bounds} below).
Using in addition \eqref{eq:2-2-est}, \eqref{eq:D-2-est}, \eqref{eq:div-2-est}, we then have, for $f\in L^2(\Om)^d$, $\vrh\in\R$, $\ph\in\E$, and $t>0$,
\beq
 &&\nn{e^{-\vrh\phi}e^{t\Delta_H}(e^{\vrh\phi}f)}{L^{q_0}(\Om)^d}\\
 &\lesssim& \nn{e^{-\vrh\phi}e^{t\Delta_H}(e^{\vrh\phi}f)}{V(\Om)}\\
 &\lesssim& \nn{D_-(e^{-\vrh\phi}e^{t\Delta_H}e^{\vrh\phi}f)}{L^2(\Om)^{d\times d}}
                          +\nn{\dv (e^{-\vrh\phi}e^{t\Delta_H}e^{\vrh\phi}f)}{L^2(\Om)}
                          +\nn{e^{-\vrh\phi}e^{t\Delta_H}(e^{\vrh\phi}f)}{L^{2}(\Om)^d}\\
 &\lesssim& (1+t^{-1/2})e^{\om_1\vrh^2t}\nn{f}{L^2(\Om)^d}.
\eeq
Hence we find for any $\del>0$ a constant $C_\del>0$ such that, for all $\vrh\in\R$, $\ph\in\E$, and $t>0$, we have
\begin{equation}\label{eq:2-q0-est}
  \nn{e^{-\vrh\phi}e^{t\Delta_H}e^{\vrh\phi}f}{L^2(\Om)^d\to L^{q_0}(\Om)^d} \le C_\del t^{-1/2}e^{\del t} e^{\om_1\vrh^2t} = C_\del t^{-\frac{d}2(\frac12-\frac1{q_0})}e^{\del t} e^{\om_1\vrh^2t}.
\end{equation}
\emph{Step 4:} We use the arguments in \cite{BK-mr} and obtain, for any $\del>0$, a constants $C_{q_0,\del},\om_{q_0}>0$ such that, for all $\vrh\in\R$, $\ph\in\E$, and $t>0$, we have
\begin{equation}\label{eq:q0-q0-est}
  \nn{e^{-\vrh\phi}e^{t\Delta_H}e^{\vrh\phi}f}{L^{q_0}(\Om)^d\to L^{q_0}(\Om)^d} 
  \le C_{q_0,\del} e^{\del t} e^{\om_{q_0}\vrh^2t}.
\end{equation}
\emph{Step 5:} We use Proposition \ref{prop:Lap-PS} for $q=q_0$ and obtain a constants $C_{q_0},\del_{q_0}>0$ such that,
for all $t>0$,
\begin{equation}\label{eq:W2-q0}
 \nn{e^{t\Delta_H}}{L^{q_0}(\Om)^d\to W^{2,q_0}(\Om)^d}\le C_{q_0} t^{-1} e^{\del_{q_0}t}.
\end{equation}
Then we use Stein interpolation between \eqref{eq:2-q0-est} and \eqref{eq:W2-q0} and obtain new constants $C_{q_0},\del_{q_0}>0$ such that, for all $\vrh\in\R$, $\ph\in\E$, and $t>0$,
\begin{equation}\label{eq:W1-q0-est}
  \nn{e^{-\vrh\phi}e^{t\Delta_H}e^{\vrh\phi}f}{L^{q_0}(\Om)^d\to W^{1,q_0}(\Om)^d} 
  \le C_{q_0,\del} t^{-1/2} e^{\del_{q_0} t} e^{\om_{q_0}\vrh^2t}.
\end{equation}
\emph{Step 6:} If $q_0<d$ we use the Sobolev embedding $W^{1,q_0}(\Om)^d\emb L^{q_1}(\Om)^d$ as in step 3, where $\frac1{q_1}=\frac1{q_0}-\frac1d$, and obtain constants $C_{q_1},\del_{q_1},\om_{q_1}>0$ such that, for all $\vrh\in\R$, $\ph\in\E$, and $t>0$, we have
\begin{equation}\label{eq:q0-q1-est}
  \nn{e^{-\vrh\phi}e^{t\Delta_H}e^{\vrh\phi}f}{L^{q_0}(\Om)^d\to L^{q_1}(\Om)^d} 
  \le C_{q_1} t^{-1/2}e^{\del_{q_1} t} e^{\om_{q_1}\vrh^2t} = C_{q_1} t^{-\frac{d}2(\frac1{q_0}-\frac1{q_1})}e^{\del_{q_1} t} e^{\om_{q_1}\vrh^2t}.
\end{equation}
Combining \eqref{eq:2-q0-est} and \eqref{eq:q0-q1-est} via the semigroup property we obtain constants $C_{q_0,q_1},\del_{q_0,q_1},\om_{q_0,q_1}>0$ such that, for all $\vrh\in\R$, $\ph\in\E$, and $t>0$, we have
\begin{equation}\label{eq:2-q1-est}
  \nn{e^{-\vrh\phi}e^{t\Delta_H}e^{\vrh\phi}f}{L^{2}(\Om)^d\to L^{q_1}(\Om)^d} \le
   C_{q_0,q_1} t^{-\frac{d}2(\frac12-\frac1{q_1})}e^{\del_{q_0,q_1} t} e^{\om_{q_0,q_1}\vrh^2t}.
\end{equation}
and can repeat steps 4--6 with $q_1$ in place of $q_0$.

If $q_0=d$ then the Sobolev embedding into $L^\infty(\Om)^d$ is not available. We interpolate between \eqref{eq:2-2-est} and \eqref{eq:q0-q0-est} to obtain \eqref{eq:q0-q0-est} for some $2<\wt{q}_0<d$ and can repeat steps 5 and 6. 

If $q_0>d$ we use the Gagliardo-Nirenberg inequality
$$
 \nn{u}{L^\infty(\Om)}\le C_{GN}\nn{u}{W^{1,q_0}(\Om)}^{d/q_0} \nn{u}{L^{q_0}(\Om)}^{1-d/q_0}
$$
instead of the Sobolev inequality and use both \eqref{eq:q0-q0-est} and \eqref{eq:W1-q0-est}. This yields
constants $C_{\infty},\del_{\infty},\om_{\infty}>0$ such that, for all $\vrh\in\R$, $\ph\in\E$, and $t>0$, we have
\eqref{eq:q0-q1-est} with $q_1=\infty$. Again, we can combine this with \eqref{eq:2-q0-est} and obtain 
constants $C_{q_0,\infty},\del_{q_0,\infty},\om_{q_0,\infty}>0$ such that, for all $\vrh\in\R$, $\ph\in\E$, and $t>0$, we have \eqref{eq:2-q1-est} for $q_1=\infty$, i.e.,
\begin{equation}\label{eq:2-infty-est}
  \nn{e^{-\vrh\phi}e^{t\Delta_H}e^{\vrh\phi}f}{L^{2}(\Om)^d\to L^{\infty}(\Om)^d} \le
   C_{q_0,\infty} t^{-\frac{d}4}e^{\del_{q_0,\infty} t} e^{\om_{q_0,\infty}\vrh^2t}.
\end{equation}
Since $e^{t\Delta_H}$ is self-adjoint, dualization of \eqref{eq:2-infty-est} yields 
\begin{equation}\label{eq:1-2-est}
  \nn{e^{-\vrh\phi}e^{t\Delta_H}e^{\vrh\phi}f}{L^{1}(\Om)^d\to L^{2}(\Om)^d} \le
   C_{q_0,\infty} t^{-\frac{d}4}e^{\del_{q_0,\infty} t} e^{\om_{q_0,\infty}\vrh^2t}.
\end{equation}
Combining \eqref{eq:2-infty-est} and \eqref{eq:1-2-est} finally yields constants $\ov{C},\ov{\del},\ov{\om}>0$ such that,
for all $\vrh\in\R$, $\ph\in\E$, and $t>0$, 
\begin{equation}\label{eq:1-infty-est}
  \nn{e^{-\vrh\phi}e^{t\Delta_H}e^{\vrh\phi}f}{L^{1}(\Om)^d\to L^{\infty}(\Om)^d} \le
   \ov{C} t^{-\frac{d}2}e^{\ov{\del} t} e^{\ov{\om}\vrh^2t}.
\end{equation}
This is well-known to imply that the operators $e^{t\Delta_H}$ have integral kernels satisfying pointwise Gaussian bounds.
As the semigroup leaves $L^2(\Om;\R^d)$ invariant, the kernels can be chosen to be $\R^{d\times d}$-valued.
\end{proof}

\begin{remark}\label{rem:kernel-bounds}\rm
In case $d=2$ one has to use the Gagliardo-Nirenberg type inequality 
$$
 \nn{u}{L^{q_0}(\Om)}\le C_{GN} \nn{u}{W^{1,2}(\Om)}^{1-2/{q_0}}\nn{u}{L^2(\Om)}^{2/{q_0}}
$$
for some $2<q_0<\infty$ since the Sobolev embedding $V(\Om)\subseteq W^{1,2}(\Om)^d\emb L^{q_0}(\Om)^d$ does not give the right $t$-exponent in \eqref{eq:2-q0-est}.  
\end{remark}

We note some consequences of Theorem~\ref{thm:Gaussian-bds}.

\begin{corollary}\label{cor:Lap-Hq}
Let $\Om\subseteq\R^d$ be a uniform $C^{2,1}$-domain. Then the semigroup $(T(t))_{t\ge0}$ generated by $\Delta_H$ extends for $q\in[1,\infty]$ to consistent analytic semigroups $(T_q(t))_{t\ge0}$ on $L^q(\Om)^d$ whose generators we denote by $\Delta_{H,q}$. On $L^q(\Om)^d$ the corresponding semigroup is strongly continuous for $q\in[1,\infty)$, on $L^\infty(\Om)^d$ the semigroup is $w^*$-continuous, and we have the duality relation $T_q(t)^*=T_{q'}(t)$ for $t\ge0$, hence also $(\Delta_{H,q})^*=\Delta_{H,q'}$. For $q\in(1,\infty)$ we have $\Delta_{H,q}=\Delta_{PS,q}$ and thus the domain description in Proposition~\ref{prop:Lap-PS}.

The semigroup operators $T_\infty(t)$, $t>0$, leave $C_0(\ov{\Om})^d$ invariant and thus induce an analytic semigroup $(T_0(t))_{t\ge0}$ in $C_0(\ov{\Om})^d$ whose generator we denote by $\Delta_{H,0}$.
\end{corollary}

\begin{proof}
The assertion on extension of $(T(t))_{t\ge0}$ to analytic semigroups on $L^q(\Om)^d$ for $q\in[1,\infty]$ is a well-known consequence of pointwise Gaussian bounds (see, e.g., \cite{Ouh:book}). Observe also that here $T_\infty(t)=T_1(t)'$, $t\ge0$, due to self-adjointness of $\Delta_H$. By consistency we have $T_q(t)=e^{t\Delta_{PS,q}}$, $t\ge0$, for $q\in(1,\infty)$ hence
$\Delta_{H,q}=\Delta_{PS,q}$ with domain given in Proposition~\ref{prop:Lap-PS} for $q\in(1,\infty)$. 

Let $f\in C_c(\ov{\Om})^d$ and $t>0$. Choose $q>\frac{d}2$. Then $f\in L^2(\Om)^d\cap L^q(\Om)^d$ and, by analyticity in $L^q(\Om)^d$ and Sobolev embedding,
$$
 T_q(t)f\in D(\Delta_{PS,q})\subseteq W^{2,q}(\Om)^d\emb C_0(\ov{\Om})^d.
$$ 
Since $C_c(\ov{\Om})^d$ is dense in $C_0(\ov{\Om})^d$ w.r.t. to $\nn{\cdot}{\infty}$ we conclude that the operators $T_\infty(t)$, $t\ge0$, leave $C_0(\ov{\Om})^d$ invariant.
\end{proof}

\begin{remark}\rm\label{rem:tilde-Lq}
Let $\Om\subseteq\R^d$ be a uniform $C^{2,1}$-domain and $q\in(1,\infty)$.
By consistency of the semigroups $(T_2(t))_{t\ge0}$ on $L^2(\Om)^d$ and $(T_q(t))_{t\ge0}$ on $L^q(\Om)^d$ we obtain a consistent analytic semigroup $(\wt{T}_q(t))_{t\ge0}$ on $\wt{L}^q(\Om)^d$, whose generator we denote by $\wt{\Delta}_{H,q}$. Then
$$
 D(\wt{\Delta}_{H,q})=\{u\in\wt{W}^{2,q}(\Om)^d:\nu\cdot u=0\ \mbox{and}\ D_-(u)\nu=0\ \mbox{on $\partial\Om$}\,\}.
$$
With respect to the duality \eqref{eq:tilde-Lq-dual} we have $(\wt{T}_q(t))'=\wt{T}_{q'}(t)$, $t\ge0$.
\end{remark}

\begin{remark}\rm\label{rem:better-del}
The exponent $\del>0$ in \eqref{eq:kernel-bds} depends on the exponents $\del_{q_0}$ in \eqref{eq:W2-q0}, i.e. on the exponential growth of the semigroups in Proposition~\ref{prop:Lap-PS}, which is not specified in \cite[Theorem 6.1]{hobus-saal}. However,
pointwise Gaussian kernel bounds imply that the spectrum of $-\Delta_{H,q}$ does not depend on $q\in[1,\infty]$ (see, e.g., \cite{K:london}) hence equals $\si(-\Delta_{H,2})\subseteq[0,\infty)$. As the growth of an analytic semigroup is detemined by the spectral bound of its generator we find, for any $q\in[1,\infty]$ and $\eps>0$, a constant $M_{\eps,q}>0$ such that
$$
  \nn{T_q(t)}{L^q(\Om)^d\to L^q(\Om)^d}\le M_{\eps,q} e^{\eps t}\quad\mbox{for all $t>0$.}
$$  
The same holds for the growth of $(\wt{T}_q(t))_{t\ge0}$ in $\wt{L}^q(\Om)^d$ for $q\in(1,\infty)$. These improved bounds can then be used to obtain, by a repetition of the proof, an arbitrarily small $\del>0$ in \eqref{eq:kernel-bds}.
\end{remark}

Our main result on the Hodge Laplacian is as follows. 

\begin{theorem}\label{thm:HLap-Hinfty}
Let $\Om\subseteq\R^d$ be a uniform $C^{2,1}$-domain, $q\in(1,\infty)$, $\theta\in(0,\frac\pi2)$, and $\del>0$. Then the operator $\del-\Delta_{H,q}$ has a bounded $H^\infty(\Si_\theta)$-functional calculus in $L^q(\Om)^d$ and $\del-\wt{\Delta}_{H,q}$ has a bounded $H^\infty(\Si_\theta)$-functional calculus in $\wt{L}^q(\Om)^d$. 

In fact, these operators even have a H\"ormander functional calculus with an estimate as in \eqref{eq:H-Ws2-cond} for $s>(d+1)|\frac12-\frac1q|$. 
\end{theorem}

\begin{proof}
Combining Theorem~\ref{thm:Gaussian-bds} and Remark~\ref{rem:better-del} we obtain a bounded $H^\infty$-calculus for
$\del-\Delta_{H,q}$ by the main result of \cite{duong-robinson}. The result on the angle of the $H^\infty$-calculus is implied
by the much stronger H\"ormander type functional calculus that $\del-\Delta_{H,q}$ enjoys 
by the results of \cite{DOS} or \cite{KUhl:spec-mult}.
\end{proof}

\begin{remark}\rm\label{rem:Hoerm}
The arguments that led to Theorem~\ref{thm:HLap-Hinfty} are very similar to those in the applications of the results of \cite{KUhl:spec-mult} to the elliptic systems in \cite{KUhl:ell-sys}. The condition on $s$ is obtained by interpolation.
\end{remark}

As  $L^q(\Om)^d$ and  $\wt{L}^q(\Om)^d$ are UMD-spaces for $q\in(1,\infty)$, we obtain the usual consequences of a bounded $H^\infty$-calculus. 

\begin{corollary}\label{cor:HLap-BIP}
Let $\Om\subseteq\R^d$ be a uniform $C^{2,1}$-domain, $q\in(1,\infty)$, and $\del>0$.
The operators $\del-\Delta_{H,q}$ in $L^q(\Om)^d$ and $\del-\wt{\Delta}_{H,q}$ in $\wt{L}^q(\Om)^d$ 
have bounded imaginary powers.
In particular, for $\alpha\in(0,1)$, we have
\begin{equation}\label{eq:frac-c-int}
 D((\del-\Delta_{H,q})^\alpha)=[L^q(\Om)^d,D(\Delta_{H,q})]_\alpha,\qquad
 D((\del-\wt{\Delta}_{H,q})^\alpha)=[\wt{L}^q(\Om)^d,D(\wt{\Delta}_{H,q})]_\alpha.
\end{equation}
Moreover, the operators $\Delta_{H,q}$ and $\wt{\Delta}_{H,q}$ have maximal $L^p$-regularity, $p\in(1,\infty)$, on finite intervals in $L^q(\Om)^d$ and $\wt{L}^q(\Om)^d$, respectively.
\end{corollary}

Invoking Proposition~\ref{prop:seeley-interpol} we can now identify the fractional domain spaces.

\begin{corollary}\label{cor:HLap-frac-dom}
Let $\Om\subseteq\R^d$ be a uniform $C^{2,1}$-domain and $q\in(1,\infty)$. Then we have 
\begin{align*}  
 [L^q(\Om)^d, D(\Delta_{H,q})]_\al
 =\left\{\begin{array}{ll}
 H^{2\al,q}(\Om)^d, &\al\in(0,\frac1{2q}), \\
 \{u\in H^{2\al,q}(\Om)^d: \nu\cdot u|_{\partial\Om}=0 \}, &\al\in (\frac1{2q},\frac12+\frac1{2q}), \\
 \{u\in H^{2\al,q}(\Om)^d: \nu\cdot u|_{\partial\Om}=0, D_-(u)\nu|_{\partial\Om}=0  
 \}, &\al\in (\frac12+\frac1{2q},1).
 \end{array}\right.
\end{align*}
\end{corollary}

For a description in case $\al\in\{\frac1{2q},1+\frac1{2q}\}$ we refer to \cite{Seeley}.

\section{The Stokes operator with Hodge boundary conditions}\label{sec:Hodge-Stokes}

\subsection{Invariance for $q=2$}
We start with the case $q=2$ and a Lipschitz domain $\Om\subseteq\R^d$. Recall that we have the Helmholtz projection $\PP_2$ corresponding to the orthogonal decomposition $L^2_\si(\Om)\oplus G^2(\Om)$ and the projection $\mathcal{P}_2$ corresponding to the orthogonal decomposition $L^2(\Om)^d=\mathcal{L}^2_\si(\Om)\oplus \mathcal{G}^2(\Om)$.

\begin{prop}\label{prop:inv-L2}
Let $\Om\subseteq\R^d$ be a Lipschitz domain. Then $L^2_\si(\Om)$ and $\mathcal{L}^2_\si(\Om)$ are invariant under the semigroup $(T(t))_{t\ge0}$ generated by $\Delta_H$ in $L^2(\Om)^d$. 
\end{prop}

\begin{proof}
We use \cite[Theorem 2.1]{Ouh:inv-convex} and thus have to check that $u\in V(\Om)$ implies $\PP_2u,\mathcal{P}_2u\in V(\Om)$ and $\Rep\a (u,u-\PP_2u)\ge0$, $\Rep\a(u,u-\mathcal{P}_2u)\ge0$. 
First we show $\PP_2u\in V(\Om)$. We have $\PP_2u\in L^2_\si(\Om)\subseteq\{v\in L^2(\Om)^d: \dv v=0, \nu\cdot v|_{\partial\Om}=0\}=\mathcal{L}^2_\si(\Om)$ and $\mathcal{P}_2u\in\mathcal{L}^2_\si(\Om)$, and it rests to prove $D_-(\PP_2u), D_-(\mathcal{P}_2u)\in L^2(\Om)^{d\times d}$. To this end write $u=v+\nabla\psi$ where
$v\in L^2_\si(\Om)$ and $\nabla\psi\in G^2(\Om)$. Then we have, distributionally,
$$
 D_-(\PP_2 u)=D_-(v)=D_-(u)-D_-(\nabla\psi)=D_-(u)\in L^2(\Om)^{d\times d}.
$$
Similarly, writing $u=\wt{v}+\nabla\wt{\psi}$ where $\wt{v}\in\mathcal{L}^2_\si(\Om)$ and $\nabla\wt{\psi}\in\mathcal{G}^2(\Om)$, we have
$$
 D_-(\mathcal{P}_2u)=D_-(\wt{v})=D_-(u)-D_-(\nabla\wt{\psi})=D_-(u)\in L^2(\Om)^{d\times d}.
$$
We conclude $\PP_2u, \mathcal{P}_2u\in V(\Om)$ and, for $w\in\{ \PP_2u, \mathcal{P}_2u\}$,
$$
 \a (u,u-w)=\frac12\int_\Om D_-(u):\ov{D_-(u-w)}+\dv u\,\ov{\dv(u-w)}\,dx
                                     =\int_\Om |\dv u|^2\,dx \ge0,
$$
which ends the proof.
\end{proof}

\begin{remark}\label{rem:P-commute}\rm
Once we have that $V(\Om)$ is invariant under $\PP_2$ and $\mathcal{P}_2$, we might just as well have argued as in \cite[Lemma 5.4]{KUhl:ell-sys} and check directly that $\PP_2$ and $\mathcal{P}_2$ commute with $-\Delta_H$, since for $u\in D(-\Delta_H)$ and $v\in V(\Om)$ we have
$$
 \sk{\PP_2(-\Delta_H)u}{v}_{L^2(\Om)^d}=\sk{-\Delta_Hu}{\PP_2v}_{L^2(\Om)^d}=\a (u,\PP_2v)
 =\frac12\int_\Om D_-(u):\ov{D_-(v)}\,dx=\a (\PP_2u,v),
$$
which means $\PP_2u\in D(-\Delta_H)$ and $(-\Delta_H)\PP_2u=\PP_2(-\Delta_H)u$. This implies that $\PP_2$ commutes with resolvents of $\Delta_H$ and thus also with the semigroup operators $T(t)$, $t\ge0$. The argument for $\mathcal{P}_2$ is the same.
\end{remark}

\subsection{Invariance for $q\neq 2$}\label{subsec:inv-Lq}
We now consider $q\in(1,\infty)$ and a uniform $C^{2,1}$-domain $\Om\subseteq\R^d$. It is no surprise that the $L^q$-theory is more subtle. However, we have invariance of certain $L^q$-spaces of solenoidal vector fields without additional assumptions and obtain some information even for the limit cases $q=1$ and $q=\infty$. We first define these spaces on more general domains. 

\begin{definition}\rm
For an arbitrary domain $\Om\subseteq\R^d$ and $q\in[1,\infty)$ we set
$$
 \check{L}^q_\si(\Om):=\ov{L^2_\si(\Om)\cap L^q(\Om)^d}^{L^q(\Om)^d}
$$
and define
$$
 \check{C}_{0,\si}(\ov{\Om}):=\ov{ L^2_\si(\Om)\cap C_0(\ov{\Om})^d}^{\nn{\cdot}{\infty}}.
$$
For a Lipschitz domain $\Om\subseteq\R^d$  and $q\in(1,\infty)$ we set
$$
 \check{\mathcal{L}}^q_\si(\Om)=\ov{\mathcal{L}^2_\si(\Om)\cap L^q(\Om)^d}^{L^q(\Om)^d}.
$$
\end{definition}

\begin{remark}\rm
Note that $C_{0,\si}(\Om):=\ov{C^\infty_{c,\si}(\Om)}^{\nn{\cdot}{\infty}}$, a space considered in the context of Dirichlet (or ``no-slip'') boundary conditions, is not suitable here, as $u=0$ on $\partial\Om$ for any $u\in C_{0,\si}(\Om)$. Recall that, for $u\in C_0(\ov{\Om})$ we only have that $u(x)\to0$ as $|x|\to\infty$ with $x\in\ov{\Om}$, see the beginning of Subsection~\ref{subsec:function-spaces}.
\end{remark}

\begin{lemma}\label{lem:Lq-solenoidal}
Let $\Om\subseteq\R^d$ be a Lipschitz domain and $q\in(1,\infty)$. Then we have
$$
 L^q_\si(\Om)\subseteq \check{L}^q_\si(\Om)\subseteq \check{\mathcal{L}}^q_\si(\Om)\subseteq\mathcal{L}^q_\si(\Om)
$$
If $L^q_\si(\Om)=\mathcal{L}^q_\si(\Om)$ then all these spaces coincide.
If $q\in I_\PP$ then $\check{L}^q_\si(\Om)=L^q_\si(\Om)$.
If $q\in[1,\frac{d}{d-1}]$ then $L^q_\si(\Om)=\mathcal{L}^q_\si(\Om)$. 
\end{lemma}

\begin{proof}
For the first inclusion observe $C^\infty_{c,\si}(\Om)\subseteq L^2_\si(\Om)\cap L^q(\Om)^d$ and recall the definition of $L^q_\si(\Om)$. For the second inclusion recall $L^2_\si(\Om)\subseteq\mathcal{L}^2_\si(\Om)$. For the third inclusion observe that, essentially by definition,
$$
 \mathcal{L}^2_\si(\Om)\cap L^q(\Om)^d=\{f\in L^2(\Om)^d\cap L^q(\Om)^d: \dv f=0, \nu\cdot f|_{\partial\Om}=0\,\}
 =L^2(\Om)^d\cap\mathcal{L}^q_\si(\Om).
$$
Now let $q\in I_\PP$ and $f\in L^2_\si(\Om)\cap L^q(\Om)^d$. Then $f=\PP_2f=\PP_q f\in L^q_\si(\Om)$, hence
$L^2_\si(\Om)\cap L^q(\Om)^d\subseteq L^q_\si(\Om)$ and the assertion follows. 

The last assertion holds by \cite[Theorem 5]{mas-bog-appr}.
\end{proof}

\begin{prop}\label{prop:inv-Lq}
Let $\Om\subseteq\R^d$ be a uniform $C^{2,1}$-domain and $q\in(1,\infty)$. 
Then we have:
\begin{enumerate}[label=$(\roman*)$]
\item\label{item:inv-Lq-1}The spaces $\check{{L}}^q_\si(\Om)$ and $\check{\mathcal{L}}^q_\si(\Om)$ are invariant under the semigroup $(T_q(t))_{t\ge0}$. Moreover, $L^1_\si(\Om)=\mathcal{L}^1_\si(\Om)$ is invariant under $(T_1(t))_{t\ge0}$ and $\check{C}_{0,\si}(\ov{\Om})$ is invariant under $(T_0(t))_{t\ge0}$.
\item\label{item:inv-Lq-2} The space $\wt{L}^q_\si(\Om)$ is invariant under the semigroup $(\wt{T}_q(t))_{t\ge0}$ and the Helmholtz projection $\wt{P}_q$ commutes with the semigroup operators. 
\item\label{item:inv-Lq-3} If $q\in (1,\frac{d}{d-1}]\cup I_\PP\cup[2,\infty)$ then $L^q_\si(\Om)$ is invariant under the semigroup $(T_q(t))_{t\ge0}$. 
\end{enumerate}
\end{prop}

\begin{proof}
\ref{item:inv-Lq-1} We start with $\check{L}^q_\si(\Om)$. Let $t>0$. It suffices to show $u:=T(t)f\in \check{L}^q_\si(\Om)$ for $f\in L^2_\si(\Om)\cap L^q(\Om)^d$. This is clear by Proposition \ref{prop:inv-L2} and boundedness of the semigroup operator in $L^q(\Om)^d$. The proof for $\check{\mathcal{L}}^q_\si(\Om)$ is along the same lines and uses invariance of $\mathcal{L}^2_\si(\Om)$. Also the proofs for $L^1_\si(\Om)=\mathcal{L}^1_\si(\Om)=\check{L}^1_\si(\Om)$ and $\check{C}_{0,\si}(\ov{\Om})$ are similar.

\ref{item:inv-Lq-2} Here we make use of the $\wt{L}^q$-theory in \cite{FKS:ArchMath} (see Theorem~\ref{thm:FKS-dec}).
Let $\ph\in C^\infty_{c,\si}(\Om)$. By Proposition~\ref{prop:inv-L2}, for any $t>0$, we have $T(t)\ph\in L^2_\si(\Om)$. 
For $q\le2$ we immediately obtain $T(t)\ph \in \wt{L}^q_\si(\Om)$. 

For $q>2$ we obtain $T(t)\ph=\PP_2T(t)\ph=\wt{P}_qT(t)\ph\in \wt{L}^q_\si(\Om)$ via Theorem~\ref{thm:FKS-dec}. Since $C^\infty_{c,\si}(\Om)$ is dense in $\wt{L}^q_\si(\Om)$ by Theorem~\ref{thm:FKS-dec}, we obtain invariance of $\wt{L}^q(\Om)$ under $(\wt{T}_q(t))_{t\ge0}$ by Remark~\ref{rem:tilde-Lq}.

Now combine duality of semigroups in Remark~\ref{rem:tilde-Lq} with the annihilator relations in Theorem~\ref{thm:FKS-dec} to obtain invariance of $\wt{G}^q(\Om)$ under $(\wt{T}_q(t))_{t\ge0}$. Hence $\wt{P}_q$ commutes with the semigroup.

\ref{item:inv-Lq-3} For $q\in I_\PP$ the assertion follows from \ref{item:inv-Lq-1} and Lemma~\ref{lem:Lq-solenoidal}. So 
let $q\in(2,\infty)$ and $t>0$. Observe
$$
 T_q(t)\big(L^2_\si(\Om\big)\cap L^q_\si(\Om))=\wt{T}_q(t)\big(\wt{L}^q_\si(\Om)\big)\subseteq\wt{L}^q_\si(\Om)\subseteq L^q_\si(\Om).
$$
and
$$
 L^q_\si(\Om)=\ov{C^\infty_{c,\si}(\Om)}^{L^q(\Om)^d}\subseteq\ov{L^2_\si(\Om)\cap L^q_\si(\Om)}^{L^q(\Om)^d}\subseteq L^q_\si(\Om)^.
$$
Then boundedness of $T_q(t)$ on $L^q(\Om)^d$ yields $T_q(t)(L^q_\si(\Om))\subseteq L^q_\si(\Om)$ as claimed.

For $q\in(1,\frac{d}{d-1}]$ we have $L^q_\si(\Om)=\mathcal{L}^q_\si(\Om)$ by Lemma~\ref{lem:Lq-solenoidal}
and invariance follows from \ref{item:inv-Lq-1}.
\end{proof}

\begin{remark}\label{rem:HS-inv}\rm
(a) We compare Proposition~\ref{prop:inv-Lq}~\ref{item:inv-Lq-1} to the corresponding result in \cite{hobus-saal}. 
With respect to invariance for a fixed $q$ it is essentially shown in \cite[Lemma 7.2]{hobus-saal} that $T_q(t)$ maps $L^q_\si(\Om)$ into $\mathcal{L}^q_\si(\Om)$. Hence invariance of $L^q_\si(\Om)$ is obtained assuming $L^q_\si(\Om)=\mathcal{L}^q_\si(\Om)$, see \cite[Assumption 2.4]{hobus-saal} and the discussion in \cite[Remark 2.6 (c)]{hobus-saal}. We see here that it would be sufficient to assume the weaker condition $L^q_\si(\Om)=\check{L}^q_\si(\Om)$,
which by Lemma~\ref{lem:Lq-solenoidal} is implied by $q\in I_\PP$. 
However, in \ref{item:inv-Lq-3} we obtain invariance of $L^q_\si(\Om)$ for $2\le q<\infty$ without additional assumptions. Notice that we needed \ref{item:inv-Lq-2} as an intermediate step.

(b) We consider it unlikely to have invariance of $L^q_\si(\Om)$ for $q\in(1,2)$ in the general case.

(c) In the following we concentrate on the spaces $L^q_\si(\Om)$ and $\check{L}^q_\si(\Om)$ although results similar to the $\check{L}^q_\si$-case hold for the spaces $\check{\mathcal{L}}^q_\si(\Om)$ as well, and by the same methods.
\end{remark}

Proposition~\ref{prop:inv-Lq} allows us to define the following Stokes operators in solenoidal $L^q$-spaces.

\begin{definition}\rm\label{def:H-Stokes}
Let $\Om\subseteq\R^d$ be a uniform $C^{2,1}$-domain.
For $q\in(1,\infty)$ we denote by $(\check{S}_q(t))_{t\ge0}:=(T_q(t)|_{\check{L}^q_\si(\Om)})_{t\ge0}$ the Hodge Stokes semigroup in $\check{L}^q_\si(\Om)$ and by $\check{A}_{H,q}$ its negative generator, the Hodge Stokes operator in $\check{L}^q_\si(\Om)$

For $q\in(1,\infty)$ we denote by $(\wt{S}_q(t))_{t\ge0}:=(\wt{T}_q(t)|_{\wt{L}^q_\si(\Om)})_{t\ge0}$ the Hodge Stokes semigroup in $\wt{L}^q_\si(\Om)$ and by $\wt{A}_{H,q}$ its negative generator, the Hodge Stokes operator in $\wt{L}^q_\si(\Om)$.

Whenever $L^q_\si(\Om)$ is invariant under $(T_q(t))_{t\ge0}$, so in particular for all $q\in[1,\frac{d}{d-1}]\cup I_\PP\cup[2,\infty)$, we denote by 
$(S_q(t))_{t\ge0}:=(T_q(t)|_{L^q_\si(\Om)})_{t\ge0}$ the Hodge Stokes semigroup in $L^q_\si(\Om)$ and by $A_{H,q}$ its negative generator, the Hodge Stokes operator in $L^q_\si(\Om)$.

Finally, we denote by $(\check{S}_0(t))_{t\ge0}:=(T_0(t)|_{\check{C}_{0,\si}(\Om)})_{t\ge0}$ the Hodge Stokes semigroup in $\check{C}_{0,\si}(\ov{\Om})$ and by $\check{A}_{H,0}$ its negative generator, the Hodge Stokes operator in $\check{C}_{0,\si}(\ov{\Om})$.
\end{definition}

An application of Lemma~\ref{lem:inv-gen} yields the following description of the respective domains. 

\begin{prop}\label{prop:dom-H-Stokes}
Let $\Om\subseteq\R^d$ be a uniform $C^{2,1}$-domain. We have the following descriptions for the domains of the Hodge Stokes operators introduced in Definition~\ref{def:H-Stokes}.

For $q\in(1,\infty)$ we have
$$
 D(\check{A}_{H,q})=\{u\in W^{2,q}(\Om)^d\cap \check{L}^q_\si(\Om) :  D_-(u)\nu=0\mbox{ on $\partial\Om$}\,\}
$$
and 
$$
 D(\wt{A}_{H,q})=\{u\in \wt{W}^{2,q}(\Om)^d\cap \wt{L}^q_\si(\Om) : D_-(u)\nu=0\mbox{ on $\partial\Om$}\,\}.
$$
Whenever $L^q_\si(\Om)$ is invariant under $(T_q(t))_{t\ge0}$, so in particular for all $q\in (1,\frac{d}{d-1}]\cup I_\PP\cup[2,\infty)$, we have
$$
 D(A_{H,q})=\{u\in W^{2,q}(\Om)^d\cap L^q_\si(\Om) : 
 D_-(u)\nu=0\mbox{ on $\partial\Om$}\,\}
$$
Any of these operators acts on its domain as the negative distributional Laplacian $-\Delta$.

Finally we have 
$$
 D({A}_{H,1})=\{u\in D(\Delta_{H,1})\cap{L}^1_\si(\Om):\Delta_{H,1}u\in{L}^1_\si(\Om)\}
$$
and ${A}_{H,1}=\Delta_{H,1}|_{D({A}_{H,1})}$, and 
$$
  D(\check{A}_{H,0})=\{u\in D(\Delta_{H,0})\cap\check{C}_{0,\si}(\Om):\Delta_{H,0}u\in\check{C}_{0,\si}(\Om)\}
$$
and $\check{A}_{H,0}=\Delta_{H,0}|_{D(\check{A}_{H,0})}$.
\end{prop}
\begin{proof}
Combine Lemma~\ref{lem:inv-gen} with Proposition~\ref{prop:inv-Lq} and with Corollary~\ref{cor:Lap-Hq}. 
Observe that, for $q\in(1,\infty)$, any $u$ in $L^q_\si(\Om)$, $\check{L}^q_\si(\Om)$, or $\wt{L}^q_\si(\Om)$ satisfies $\nu\cdot u=0$ on $\partial\Om$.
\end{proof}

Concerning duality we have the following.

\begin{prop}\label{prop:dual-HodgeStokes}
Let $\Om\subseteq\R^d$ be a uniform $C^{2,1}$-domain.
The operator $A_{H,2}$ is self-adjoint in $L^2_\si(\Om)$ and $A_{H,2}\ge0$. For $q\in I_\PP$ we have $S_q(t)'=S_{q'}(t)$, $t\ge0$, and for $q\in(1,\infty)$ we have $\wt{S}_q(t)'=\wt{S}_{q'}(t)$, $t\ge0$.
\end{prop}
\begin{proof}
Use self-adjointness of $(T(t))_{t\ge0}$ in $L^2(\Om)^d$ and the fact that respective Helmholtz projections commute with the semigroup operators generated by the Hodge Laplacians.
\end{proof}

By restricting the functional calculi in Theorem~\ref{thm:HLap-Hinfty} to invariant subspaces we obtain our main result on Hodge Stokes operators. 

\begin{theorem}\label{thm:HS-Hinfty}
Let $\Om\subseteq\R^d$ be a uniform $C^{2,1}$-domain, $q\in(1,\infty)$, $\del>0$, and $\theta\in(0,\frac{\pi}2)$. Then
$\del+\check{A}_{H,q}$ has a bounded $H^\infty(\Si_\theta)$-calculus in $\check{L}^q_\si(\Om)$ and 
$\del+\wt{A}_{H,q}$ has a bounded $H^\infty(\Si_\theta)$-calculus in $\wt{L}^q_\si(\Om)$. 

If, in addition, $L^q_\si(\Om)$ is invariant under $(T_q(t))_{t\ge0}$, so in particular if $q\in (1,\frac{d}{d-1}]\cup I_\PP\cup [2,\infty)$, then $\del+A_{H,q}$ has a bounded $H^\infty(\Si_\theta)$-calculus in $L^q_\si(\Om)$.

In fact, these operators have a H\"ormander functional calculus with an estimate as in \eqref{eq:H-Ws2-cond} for $s>(d+1)|\frac12-\frac1q|$. 
\end{theorem}

\begin{proof}
Invariance of a closed subspace under the semigroup implies invariance under the resolvents of the generator, at least on the connected component of the resolvent set that contains a right half plane. This in turn implies invariance of the closed subspace under the operators of the $H^\infty$-calculus, see the definition in Subsection~\ref{subsec:Hinfty}.

Actually, also the operators in the H\"ormander functional calculus leave invariant a subspace that is left invariant under the semigroup. Hence the operators in Theorem~\ref{thm:HS-Hinfty} even have a H\"ormander functional calculus in the respective spaces of solenoidal vector fields.
\end{proof}

\begin{corollary}
Let $\Om\subseteq\R^d$ be a uniform $C^{2,1}$-domain, $q\in(1,\infty)$, and $\del>0$.
The operators $\del+\check{A}_{H,q}$ in $\check{L}^q_\si(\Om)$ and 
$\del+\wt{A}_{H,q}$ in $\wt{L}^q_\si(\Om)$ have bounded imaginary powers.
In particular, for $\alpha\in(0,1)$, we have
$$
 D((\del+\check{A}_{H,q})^\alpha)=[\check{L}^q_\si(\Om),D(\check{A}_{H,q})]_\alpha,\qquad
 D((\del+\wt{A}_{H,q})^\alpha)=[\wt{L}^q_\si(\Om),D(\wt{A}_{H,q})]_\alpha.
$$
Moreover, the operators $\check{A}_{H,q}$ and $\wt{A}_{H,q}$ have maximal $L^p$-regularity, $p\in(1,\infty)$, on finite intervals in $\check{L}^q_\si(\Om)$ and $\wt{L}^q_\si(\Om)$, respectively.
 
If, in addition, $L^q_\si(\Om)$ is invariant under $(T_q(t))_{t\ge0}$, in particular if $q\in (1,\frac{d}{d-1}]\cup I_\PP\cup [2,\infty)$, then $\del+A_{H,q}$ has the respective properties in $L^q_\si(\Om)$.
\end{corollary}

Combining Corollary~\ref{cor:HLap-frac-dom} with Corollary~\ref{cor:inv-fract} we obtain the following representations for the fractional domain spaces of the Hodge Stokes operator in $L^q_\si(\Om)$.

\begin{corollary}\label{cor:HStokes-frac-dom}
Let $\Om\subseteq\R^d$ be a uniform $C^{2,1}$-domain and $q\in(1,\frac{d}{d-1}]\cup I_\PP\cup[2,\infty)$ (or assume more generally that $L^q_\si(\Om)$ invariant under $(T_q(t))$. Then we have 
\begin{align}\label{eq:HStokes-frac-dom}
 [L^q_\si(\Om)^d, D(A_{H,q})]_\al
 =\left\{\begin{array}{ll}
 H^{2\al,q}(\Om)^d\cap L^q_\si(\Om) &, \ \al\in(0,\frac1{2q}), \\
 H^{2\al,q}(\Om)^d\cap L^q_\si(\Om) &, \ \al\in (\frac1{2q},\frac12+\frac1{2q}), \\
 \{u\in H^{2\al,q}(\Om)^d\cap L^q_\si(\Om): D_-(u)\nu|_{\partial\Om}=0 \} &, \ \al\in (\frac12+\frac1{2q},1).
 \end{array}\right.
\end{align}
\end{corollary}

For information on the limit cases $\al\in\{\frac1{2q},1+\frac1{2q}\}$ we refer again to \cite{Seeley}.

\section{Robin Stokes as perturbations of Hodge Stokes}\label{sec:pert-HBC}

In this section we shall perturb the Hodge boundary conditions on a uniform $C^{2,1}$-domain $\Om\subseteq\R^d$.  
This can be done in the spaces $L^q(\Om)^d$ but even for $q\in I_\PP$ the perturbed semigroup will not leave $L^q_\si(\Om)$ invariant. Hence we shall perturb the Hodge Stokes operator in $L^q_\si(\Om)$ directly. Perturbation of boundary conditions is a subtle business. In order to have precise domain descriptions we need information on the resolvent problem for the Hodge Stokes operator with inhomogeneous boundary conditions. Similar to what has been done in \cite{hobus-saal}, we shall get them from the estimates on the resolvent problem for the Hodge Laplacian with inhomogeneous boundary conditions. However, we can dispense with \cite[Assumption 2.4]{hobus-saal} which may be phrased as $L^q_\si(\Om)=\mathcal{L}^q(\Omega)$ and which has been crucial for the results in \cite{hobus-saal}, see also Remark~\ref{rem:HS-inv}.

\subsection{Estimates for resolvent problems}
We start by recalling \cite[Theorem 6.1]{hobus-saal}: Let $\Om\subseteq\R^d$ be a uniform $C^{2,1}$-domain, $q\in(1,\infty)$, $\theta\in(0,\pi)$ and $\del>0$: For any $f\in L^q(\Om)^d$, $g\in W^{1,q}(\Om)^d$, and $\la\in\del+\Si_\theta$ the problem
\begin{equation}\label{eq:Lap-u-f-g}
\begin{cases}
\ \la u -\Delta u & = \ f \quad\mbox{in $\Om$}, \\
\ D_-(u)\nu & = \ g_{\tan}\quad \mbox{on $\partial\Om$}, \\
\ \nu\cdot u & = \ 0 \quad \mbox{on $\partial\Om$},
\end{cases}
\end{equation}
has a unique solution $u\in W^{2,q}(\Om)^d$, and we have the estimate
\begin{align}\label{eq:Lap-u-f-g-est}
 \nn{\la u,\la^{1/2} \nabla u,\nabla^2 u}{L^q(\Om)}\lesssim \nn{f,\la^{1/2} g,\nabla g}{L^q(\Om)}.
\end{align}
In the following, we shall denote the unique solution of \eqref{eq:Lap-u-f-g} by 
\begin{align}
 u=R_\la f+S_\la g\quad \mbox{where}\quad R_\la f=(\la-\Delta_{H,q})^{-1}f. 
\end{align}
Notice that, if $f\in L^q_\si(\Om)$ and $L^q_\si(\Om)$ is invariant under $(T_q(t))$, then $R_\la f=(\la+A_{H,q})^{-1}f$.

We first state a lemma on invariance and regularity of decompositions.

\begin{lemma}\label{lem:inv-dec-res}
Let $\Om\subseteq\R^d$ be a uniform $C^{2,1}$-domain and $q\in (1,\infty)$. Then we have the following.
\begin{enumerate}[label=$(\roman*)$]
\item\label{item:inv-dec-1} If $u\in W^{2,q}(\Om)^d\cap L^q_\si(\Om)$ with $D_-(u)\nu=0$ on $\partial\Om$ then $\Delta u\in L^q_\si(\Om)$.
\item\label{item:inv-dec-2} If $u\in W^{2,q}(\Om)^d\cap G^q(\Om)$ then $\Delta u\in G^q(\Om)$. 
\item\label{item:inv-dec-3} Let $\theta\in(0,\pi)$, $\del>0$, $\la\in\del+\Si_\theta$, $f\in L^q(\Om)^d$, $g\in W^{1,q}(\Om)^d$, and denote by $u\in W^{2,q}(\Om)^d$ the unique solution of \eqref{eq:Lap-u-f-g}. Suppose that $u=u_0+\nabla\psi$ with $u_0\in L^q_\si(\Om)$ and $\nabla\psi\in G^q(\Om)$. Then $\nabla\psi\in D(\Delta_{H,q})$ and 
$u_0\in W^{2,q}(\Om)^d\cap L^q_\si(\Om)$ with $D_-(u_0)\nu=g_{\tan}$ on $\partial\Om$. 
\end{enumerate}
\end{lemma}

\begin{proof}
\ref{item:inv-dec-1}: For $v\in G^{q'}(\Om)$ we have $D_-(v)=0$ in $\Om$. Hence 
$$
 \int_\Om (-\Delta u)\cdot\ov{v}\,dx 
 =0
$$
by Lemma~\ref{lem:calc-proof3.4}~\ref{item:calc-pf3.4-i}. We conclude that $\Delta u\in G^{q'}(\Om)^\perp=L^q_\si(\Om)$.

\ref{item:inv-dec-2}: Let $u=\nabla\psi\in W^{2,q}(\Om)^d$. Then $\Delta u=\Delta\nabla\psi=\nabla\Delta\psi\in G^q(\Om)$.

\ref{item:inv-dec-3}: We have $D_-(\nabla\psi)=0$ and $D_-(\nabla\psi)\nu=0$ on $\partial\Om$. Hence
\begin{align}\label{eq:D-u0}
 D_-(u_0)=D_-(u)\in W^{1,q}(\Om)^{d\times d}\quad\mbox{and}\quad D_-(u_0)\nu=g_{\tan}\ \mbox{on $\partial\Om$.}
\end{align}
Further we have 
\begin{align}\label{eq:Delta-psi}
 \Delta\psi=\mbox{div}\nabla\psi=\mbox{div}\,u\in W^{1,q}(\Om)\quad\mbox{and}\quad
 \Delta\nabla\psi=\nabla\Delta\psi=\nabla\mbox{div}\,u\in G^q(\Om)\subseteq L^q(\Om)^q,
\end{align}
which implies $\Delta u_0=\Delta u-\Delta\nabla\psi\in L^q(\Om)^d$. 
Finally, $\nu\cdot u_0=0$ on $\partial\Om$ implies $\nu\cdot\nabla\psi=\nu\cdot(u-u_0)=0$ on $\partial\Om$. It now suffices to show $\nabla\psi\in W^{2,q}(\Om)^d$. We have $\Delta_{H,q'}=(\Delta_{H,q})^*$ by Corollary~\ref{cor:Lap-Hq}, and for $v\in D(\Delta_{H,q'})$ we have, by Lemma~\ref{lem:calc-proof3.4}~\ref{item:calc-pf3.4-ii},
\begin{align*}
\int_\Om \nabla\psi\cdot\ov{(\ov{\la}-\Delta)v}\,dx &=\ \la\int_\Om \nabla\psi\cdot\ov{v}\,dx
+\frac12\int_\Om D_-(\nabla\psi):\ov{D_-(v)}\,dx+\int_\Om\mbox{div}\nabla\psi\, \ov{\mbox{div}v}\,dx\\ &
=\ \int_\Om (\la \nabla\psi-\Delta\nabla\psi)\cdot\ov{v}\,dx.
\end{align*}
Hence, $\nabla\psi\in D((\ov{\la}-\Delta_{H,q'})^*)=D(\la-\Delta_{H,q})\subseteq W^{2,q}(\Om)^d$, which then implies also
$u_0=u-\nabla\psi\in W^{2,q}(\Om)^d$.
\end{proof}

We now can formulate our result on the Stokes resolvent problem. Besides invariance of $L^q_\si(\Om)$ under $(T_q(t))$ we assume the following variant of the Helmholtz decomposition.

\begin{assumption}\label{ass:HH-dec}
There exists a closed subspace $\wh{G}^q(\Om)\subseteq G^q(\Om)$ such that 
$$
 L^q(\Om)^d=L^q_\si(\Om)\oplus \wh{G}^q(\Om)
$$
as a topological sum. We denote by $\wh{P}_q$ the corresponding bounded projection in $L^q(\Om)^d$ onto $L^q_\si(\Om)$ with kernel $\wh{G}^q(\Om)$ and let $\wh{Q}_q:=I-\wh{P}_q$.
\end{assumption}

The following is our result on the Hodge Stokes resolvent system.

\begin{theorem}\label{thm:Stokes-res}
Let $\Om\subseteq\R^d$ be a uniform $C^{2,1}$-domain. Let $q\in (1,\infty)$ be such that $L^q_\si(\Om)$ is invariant under
$(T_q(t))$ and such that Assumption~\ref{ass:HH-dec} holds. Let $\theta\in (0,\pi)$, $\del>0$ and $\la\in\del+\Si_\theta$. For any $f\in L^q_\si(\Om)$ and $g\in W^{1,q}(\Om)^d$ there exists a unique solution $(u,\nabla p)\in\big(W^{2,q}(\Om)^d\cap L^q_\si(\Om)\big)\times \wh{G}^q(\Om)$ of the problem
\begin{equation}\label{eq:Stokes-u-f-g}
\begin{cases}
\ \la u -\Delta u +\nabla p& = \ f \quad\mbox{in $\Om$}, \\
\ D_-(u)\nu & = \ g_{\tan}\quad \mbox{on $\partial\Om$}, \\
\ \nu\cdot u & = \ 0 \quad \mbox{on $\partial\Om$},
\end{cases}
\end{equation}
and we have the estimate
\begin{align}\label{eq:est-HStokes-res}
  \nn{\la u,\la^{1/2} \nabla u,\nabla^2 u,\nabla p}{L^q(\Om)}\lesssim \nn{f,\la^{1/2} g,\nabla g}{L^q(\Om)}.
\end{align}
Moreover, we can represent the solution $(u,\nabla p)$ as 
\begin{align}\label{eq:Stokes-res-rep-u}
  u & =(\la+A_{H,q})^{-1}f+\wh{P}_q S_\la g- (\la+A_{H,q})^{-1}\wh{P}_q\nabla\dive S_\la g,  \\
  \label{eq:Stokes-res-rep-p}
  \nabla p & =\wh{Q}_q(\la S_\la g-\nabla\dive S_\la g).
\end{align}
\end{theorem}

\begin{proof}
\emph{Step 1}: We show uniqueness. So let $(u,\nabla p)\in 
\big((W^{2,q}(\Om)^d\cap L^q_\si(\Om)\big)\times \wh{G}^q(\Om)$ solve \eqref{eq:Stokes-u-f-g} with $f=0$ and $g=0$.
By Lemma~\ref{lem:inv-dec-res}\ref{item:inv-dec-1} we have $\la u-\Delta u\in L^q_\si(\Om)$, hence $\nabla p\in L^q_\si(\Om)\cap\wh{G}^q(\Om)=\{0\}$. We conclude $u=-(\la-\Delta_{H,q})^{-1}\nabla p=0$.

\emph{Step 2}: The case $g=0$. Since $L^q_\si(\Om)$ is invariant under $(T_q(t))$ it is also invariant under $(\la-\Delta_{H,q})^{-1}$ for $\la\in\del+\Si_\theta$. Hence the case $g=0$ is clear with $\nabla p=0$ and 
$u=(\la-\Delta_{H,q})^{-1}f=(\la+A_{H,q})^{-1}f$, and we get \eqref{eq:est-HStokes-res} from \eqref{eq:Lap-u-f-g-est}.

\emph{Step 3}: The case $f=0$. Let $g\in W^{1,q}(\Om)^d$. Denote by $\wt{u}=S_\la g$ the solution of \eqref{eq:Lap-u-f-g} with $f=0$ and put
$u_0:=\wh{P}_q\wt{u}=\wh{P}_q S_\la g$ and $\nabla\psi=\wh{Q}_q\wt{u}=\wh{Q}_q S_\la g$. By Lemma~\ref{lem:inv-dec-res}\ref{item:inv-dec-3} we have $u_0,\nabla\psi\in W^{2,q}(\Om)^d$, and $u_0$ solves
\begin{equation}\label{eq:Lap-u0-g}
\begin{cases}
\ \la u_0 -\Delta u_0 & = \ -\la\nabla\psi+\Delta\nabla\psi \quad\mbox{in $\Om$}, \\
\ D_-(u_0)\nu & = \ g_{\tan}\quad \mbox{on $\partial\Om$}, \\
\ \nu\cdot u_0 & = \ 0 \quad \mbox{on $\partial\Om$},
\end{cases}
\end{equation}
where we recall
$$
 \Delta\nabla\psi=\nabla\Delta\psi=\nabla\dive(\nabla\psi)=\nabla\dive\wt{u}=\nabla\dive S_\la g.
$$
This term on the right hand side of the first line of \eqref{eq:Lap-u0-g} might not yet be in $\wh{G}^q(\Om)$. Hence we solve
\begin{equation}\label{eq:Lap-u1-0}
\begin{cases}
\ \la u_1 -\Delta u_1 & = \ -\wh{P}_q\nabla\dive S_\la g \quad\mbox{in $\Om$}, \\
\ D_-(u_1)\nu & = \ 0\quad \mbox{on $\partial\Om$}, \\
\ \nu\cdot u_1 & = \ 0 \quad \mbox{on $\partial\Om$},
\end{cases}
\end{equation}
where 
$$
 u_1=-(\la-\Delta_{H,q})^{-1}\wh{P}_q\nabla\dive S_\la g=-(\la+A_{H,q})^{-1}\wh{P}_q\nabla\dive S_\la g\in W^{2,q}(\Om)^d\cap L^q_\si(\Om)
$$ 
by the invariance assumption.

For $u:=u_0+u_1\in W^{2,q}(\Om)^d\cap L^q_\si(\Om)$ we then have
$$
 \la u-\Delta u=-\la\nabla\psi+\nabla\dive\wt{u}-\wt{P}_q\nabla\dive\wt{u}=\wt{Q}_q\big(-\la\wt{u}+\nabla\dive\wt{u}\big)
$$
and $u$ satisfies the boundary conditions $D_-(u)\nu=g_{\tan}$ and $\nu\cdot u=0$ on $\partial\Om$.
Letting 
\begin{align}\label{eq:def-nabla-p}
 \nabla p:=\wh{Q}_q\big(\la\wt{u}-\nabla\dive\wt{u}\big)=\wh{Q}_q\big(\la S_\la g-\nabla\dive S_\la g\big)
 \in\wh{G}^q(\Om)
\end{align}
we hence have a solution $(u,\nabla p)\in\big(W^{2,q}(\Om)^d\cap L^q_\si(\Om)\big)\times\wh{G}^q(\Om)$ of \eqref{eq:Stokes-u-f-g} for $f=0$ with the representation \eqref{eq:Stokes-res-rep-u} and \eqref{eq:Stokes-res-rep-p}. 

It rests to show \eqref{eq:est-HStokes-res}. Applying \eqref{eq:Lap-u-f-g} to \eqref{eq:Lap-u0-g} we get
$$
 \nn{\la u_0,\la^{1/2}u_0,\nabla^2 u_0}{L^q(\Om)}\lesssim
 \nn{\la\wh{Q}_qS_\la g,\nabla\dive S_\la g,\la^{1/2}g,\nabla g}{L^q(\Om)},
$$
and, by \eqref{eq:Lap-u-f-g} again,
$$
 \nn{\la\wh{Q}_qS_\la g,\nabla\dive S_\la g}{L^q(\Om)}\lesssim\nn{\la^{1/2}g,\nabla g}{L^q(\Om)}.
$$
Applying \eqref{eq:Lap-u-f-g} to \eqref{eq:Lap-u1-0} we get
$$
 \nn{\la u_1,\la^{1/2}u_1,\nabla^2 u_1}{L^q(\Om)}\lesssim
 \nn{\wh{P}_q\nabla\dive S_\la g}{L^q(\Om)}\lesssim\nn{\la^{1/2}g,\nabla g}{L^q(\Om)}.
$$
Finally, we apply \eqref{eq:Lap-u-f-g} to \eqref{eq:def-nabla-p} and get
$$ 
 \nn{\nabla p}{L^q(\Om)}\lesssim\nn{\la S_\la g,\nabla\dive S_\la g}\lesssim\nn{\la^{1/2} g,\nabla g}{L^q(\Om)},
$$
which finishes the proof of \eqref{eq:est-HStokes-res}.
\end{proof} 

\begin{remark}\rm
(a) Notice that Assumption~\ref{ass:HH-dec} holds for $q\in I_\PP$ with $\wh{G}^q(\Om)=G^q(\Om)$ and $\wh{P}_q=P_q$ 
and $\wh{Q}_q=I-P_q$. By Proposition~\ref{prop:inv-Lq}\ref{item:inv-Lq-3}, $q\in I_\PP$ also implies invariance $L^q_\si(\Om)$ under $(T_q(t))$. For $q\in I_\PP$ we have 
$$
 \nabla p=\la S_\la g-\nabla\dive S_\la g
$$
in \eqref{eq:Stokes-res-rep-p} and the term $\wh{P}_q \nabla\dive S_\la g$ in \eqref{eq:Stokes-res-rep-u} vanishes. In the proof we then simply have $u_1=0$.

(b) As mentioned above the estimates for the inhomogeneous resolvent system in 
\cite[Theorem 3.3]{hobus-saal} have been shown under \cite[Assumption 2.4]{hobus-saal}. By Remark~\ref{rem:disc-As2.4}
this assumption is equivalent to $L^q_\si(\Om)=\mathcal{L}^q(\Omega)$, so it is clearly stronger than 
invariance of $L^q_\si(\Om)$ under $(T_q(t))$, see Subsection~\ref{subsec:inv-Lq}.

(c) The case $f\in L^q_\si(\Om)$ is sufficient for our purposes. Under the same assumptions one can obtain a version of Theorem~\ref{thm:Stokes-res} for general $f\in L^q(\Om)^d$. All one has to do is to replace $f$ in \eqref{eq:Stokes-res-rep-u} by $\wh{P}_qf$ and add the term $\wh{Q}_qf$ to the representation of $\nabla p$ in \eqref{eq:Stokes-res-rep-p}.
\end{remark}

\subsection{The Robin Stokes operator in $L^q_\si(\Om)$}
For $q\in(1,\infty)$ satisfying the assumptions of Theorem~\ref{thm:Stokes-res}  and $B\in C^{0,1}(\partial\Omega)^{d\times d}$ we can now define the Robin Stokes operator $A_{B,q}$ by
$$
 A_{B,q}u:= - \wh{P}_q\Delta u,\qquad u\in D(A_{B,q}),
$$
with
$$
 D(A_{B,q}):=\{u\in W^{2,q}(\Om)^d\cap L^q_\si(\Om): D_-(u)\nu=[Bu]_{\tan}\ \mbox{on}\ \partial\Om\}.
$$
The following is our main result on Robin Stokes operators in $L^q_\si(\Om)$-spaces.

\begin{theorem}\label{thm:Robin-Stokes-Lip}
Let $\Om\subseteq\R^d$ be a uniform $C^{2,1}$-domain and let $B\in C^{0,1}(\partial\Om)^{d\times d}$. 
Let $q\in (1,\infty)$ be such that $L^q_\si(\Om)$ is invariant under $(T_q(t))$ and such that Assumption~\ref{ass:HH-dec} holds. For $\theta\in(0,\frac{\pi}2)$ there exists $\del_0>0$ such that the operator $\delta+A_{B,q}$ has a bounded $H^\infty(\Si_\theta)$-calculus in $L^q_\si(\Om)$.
\end{theorem}

\begin{proof}
We extend $B$ to a Lipschitz function on $\ov{\Om}$ with $\nn{B,\nabla B}{L^\infty(\Om)}\lesssim\nn{B,\nabla B}{L^\infty(\partial\Om)}$. We fix $\del>0$.
For $f\in L^q_\si(\Om)$ and $\la\in\del+\Si_\si$ with $\theta+\frac\pi2<\si<\pi$, we study the resolvent problem
\begin{equation}\label{eq:Stokes-u-f-Bu}
\begin{cases}
\ \la u -\Delta u +\nabla p& = \ f \quad\mbox{in $\Om$}, \\
\ D_-(u)\nu & = \ [Bu]_{\tan}\quad \mbox{on $\partial\Om$}, \\
\ \nu\cdot u & = \ 0 \quad \mbox{on $\partial\Om$}.
\end{cases}
\end{equation}
via Theorem~\ref{thm:Stokes-res} and \cite[Lemma 7.10]{KuW:levico}. For $u\in W^{2,q}(\Om)^d\cap L^q_\si(\Om)$ we have
$$
 \nn{\la^{1/2} Bu,\nabla Bu}{L^q}\lesssim \nn{B}{\infty}\nn{\la^{1/2}u,\nabla u}{L^q}+\nn{\nabla B}{L^\infty}\nn{u}{L^q}
 \lesssim \la^{-1/2}\nn{\la u,\la^{1/2}\nabla u}{L^q}.
$$  
By \cite[Lemma 7.10]{KuW:levico} we infer that for $\la\in\del+\Si_\si$ with $|\la|$ sufficiently large, the problem \eqref{eq:Stokes-u-f-Bu} has a unique solution with the estimate
$$
 \nn{\la u,\la^{1/2}\nabla u,\nabla^2 u,\nabla p}{L^q(\Om)}\lesssim\nn{f}{L^q(\Om)}.
$$
We conclude that, for $\del_0>\del$ sufficiently large, $\del_0+A_{B,q}$ is sectorial in $L^q_\si(\Om)$ and $\la\in\rho(A_{B,q})$ for $\la\in\del_0+\Si_\si$ with 
$$
 (\la+A_{B,q})^{-1}f=(\la+A_{H,q})^{-1}f+\wh{P}_qS_\la B(\la+A_{B,q})^{-1}f -(\la+A_{H,q})^{-1}\wh{P}_q\nabla\dive S_\la B(\la+A_{B,q})^{-1}f
$$
and the estimates 
$$
 \nn{\la(\la+A_{B,q})^{-1}f,\la^{1/2}\nabla(\la+A_{B,q})^{-1}f,\nabla^2(\la+A_{B,q})^{-1}f}{L^q(\Om)}\lesssim\nn{f}{L^q(\Om)}.
$$
Since we then have
\begin{align*}
& \nn{\la\wh{P}_qS_\la B(\la+A_{B,q})^{-1}f ,\la(\la+A_{H,q})^{-1}\wh{P}_q\nabla\dive S_\la B(\la+A_{B,q})^{-1}f}{L^q(\Om)} \\
\lesssim& \nn{\la^{1/2} B(\la+A_{B,q})^{-1}f,\nabla B(\la+A_{B,q})^{-1}f}{L^q(\Om)} \\
\lesssim& \la^{-1/2}\nn{\la(\la+A_{B,q})^{-1}f,\la^{1/2}\nabla (\la+A_{B,q})^{-1}f}{L^q(\Om)} \\
\lesssim& \la^{-1/2} \nn{f}{L^q(\Om)},
\end{align*}
we can see directly that the contour integral over the perturbative term yields a bounded operator in $L^q_\si(\Om)$, see \eqref{eq:Hinfty-contour}. Since $\del_0+A_{H,q}$ has a bounded $H^\infty(\Si_\theta)$-calculus, we conclude that also $\del_0+A_{B,q}$ has a bounded $H^\infty$-calculus. A similar argument has been used in \cite{AHS}.
\end{proof}

\begin{corollary}\label{cor:Robin-Stokes}
Under the assumptions of Theorem~\ref{thm:Robin-Stokes-Lip} and for $\del_0>0$ large enough, the operator 
$\del_0+A_{B,q}$ has bounded imaginary powers. In particular, for $\alpha\in(0,1)$, we have
$$
 D((\del_0+{A}_{B,q})^\alpha)=[{L}^q_\si(\Om),D({A}_{B,q})]_\alpha
$$
and 
\begin{align}\label{eq:RStokes-frac-dom}
 [L^q_\si(\Om)^d, D(A_{B,q})]_\al
 =\left\{\begin{array}{ll}
 H^{2\al,q}(\Om)^d\cap L^q_\si(\Om), &\!\! \al\in(0,\frac1{2q}), \\
 H^{2\al,q}(\Om)^d\cap L^q_\si(\Om), &\!\! \al\in (\frac1{2q},\frac12+\frac1{2q}), \\
 \{u\in H^{2\al,q}(\Om)^d\cap L^q_\si(\Om): D_-(u)\nu|_{\partial\Om}=Bu \}, &\!\! \al\in (\frac12+\frac1{2q},1).
 \end{array}\right.
\end{align}
Moreover, in ${L}^q_\si(\Om)$ the operator ${A}_{B,q}$ has maximal $L^p$-regularity on finite intervals, $p\in(1,\infty)$. 
\end{corollary}

The assertions are immediate, except for \eqref{eq:RStokes-frac-dom}. For this we shall need a result for the corresponding Robin Laplacian $\Delta_{B,q}$, given by
$$
 \Delta_{B,q}u:=\Delta u,\qquad D(\Delta_{B,q}),
$$
with
$$
  D(\Delta_{B,q}):=\{u\in W^{2,q}(\Om)^d: \nu\cdot u|_{\partial\Om}=0, D_-(u)\nu|_{\partial\Om}=Bu\},
$$
which we present next.

\begin{prop}\label{prop:Robin-Lap}
Let $\Om\subseteq\R^d$ be a uniform $C^{2,1}$-domain, $B\in C^{0,1}(\partial\Om)^{d\times d}$ and $q\in(1,\infty)$.
 For $\theta\in(0,\frac{\pi}2)$ there exists $\del>0$ such that the operator $\delta-\Delta_{B,q}$ has a bounded $H^\infty(\Si_\theta)$-calculus in $L^q_\si(\Om)$. Moreover, we have 
 \begin{align}\label{eq:RLap-frac-dom}
 [L^q_\si(\Om)^d, D(\Delta_{B,q})]_\al
 =\left\{\begin{array}{ll}
 H^{2\al,q}(\Om)^d, &\!\! \al\in(0,\frac1{2q}), \\
 \{u\in H^{2\al,q}(\Om)^d: \nu\cdot u|_{\partial\Om}=0\}, &\!\! \al\in (\frac1{2q},\frac12+\frac1{2q}), \\
 \{u\in H^{2\al,q}(\Om)^d:\nu\cdot u|_{\partial\Om}=0, D_-(u)\nu|_{\partial\Om}=Bu \}, &\!\! \al\in (\frac12+\frac1{2q},1).
 \end{array}\right.
\end{align}
\end{prop}

\begin{proof}[Proof of Corollary~\ref{cor:Robin-Stokes}]
The proof is similar to the proof of Theorem~\ref{thm:Robin-Stokes-Lip} but in fact simpler, as instead of using Theorem~\ref{thm:Stokes-res} we can directly rely on the resolvent system \eqref{eq:Lap-u-f-g} and the estimate \eqref{eq:Lap-u-f-g-est}. This yields a similar resolvent estimate for the Robin Laplacian. By Seeley's result (\cite{Seeley}) again, we obtain \eqref{eq:RLap-frac-dom}.
\end{proof}

\begin{proof}[Proof of Corollary~\ref{cor:Robin-Stokes}]
We show \eqref{eq:RStokes-frac-dom}. We can get ``$\subseteq$'' by $L^q_\si(\Om)\subseteq L^q(\Om)^d$ and $D(A_{B,q})\subseteq D(\Delta_{B,q})$, $D(A_{B,q})\subseteq L^q_\si(\Om)$. Equality holds by an argument which we borrow from \cite{Giga}. We fix $\mu>\del$ and define $P_B:=\iota_q(\mu+A_{B,q})^{-1}P_q(\mu-\Delta_{B,q})$ which is a projection in $D(\Delta_{B,q})$ onto $D(A_{B,q})$. Here $\iota_q$ denotes the embedding $L^q_\si(\Om)\to L^q(\Om)$. The operator $P_B$ has a bounded extension $\wt{P}_B$ to projection in $L^q(\Om)^d$ onto $L^q_\si(\Om)$, since the dual operator $P_B^*=(\mu-\Delta_{B^*,q'})\iota_{q'}(\mu+A_{B^*,q'})^{-1}P_{q'}$ is bounded in $L^{q'}(\Om)^d$. The latter holds by
$$
 \nn{P_B^*g}{L^{q'}}\lesssim\nn{(\mu+A_{B^*,q'})^{-1}g}{W^{2,q'}}\lesssim \nn{g}{L^{q'}},
$$
where we used the estimate \eqref{eq:Lap-u-f-g-est}, but for the Robin Laplacian.
\end{proof}

\begin{remark}\rm
Theorem~\ref{thm:Robin-Stokes-Lip} and Corollary~\ref{cor:Robin-Stokes} cover Stokes operators with Navier boundary conditions as in \eqref{eq:BC-alpha-beta} if we take $B$ as specified in \eqref{eq:Navier-Robin-BC}.
\end{remark}

\subsection{The Robin Stokes operator in $\wt{L}^q_\si(\Om)$}
Let $\Om\subseteq\R^d$ be a uniform $C^{2,1}$-domain.
We have analogs of the results in the previous subsection in $\wt{L}^q_\si(\Om)$ for all $q\in (1,\infty)$. We only state the results and omit the detailed arguments but the starting point is again the system \eqref{eq:Lap-u-f-g}. From \cite[Theorem 6.1]{hobus-saal} we infer estimates
\begin{align}
\nn{\la u,\la^{1/2}\nabla u,\nabla^2u}{\wt{L}^q(\Om)}\lesssim \nn{f,\la^{1/2}g,\nabla g}{\wt{L}^q(\Om)}.
\end{align}
We can then procede as before and obtain the following.
 
 \begin{theorem}\label{thm:Robin-Stokes-Lip-tilde}
Let $\Om\subseteq\R^d$ be a uniform $C^{2,1}$-domain and let $B\in C^{0,1}(\partial\Om)^{d\times d}$ and $q\in (1,\infty)$. For $\theta\in(0,\frac{\pi}2)$ there exists $\del_0>0$ such that the operator $\delta+A_{B,q}$ has a bounded $H^\infty(\Si_\theta)$-calculus.
\end{theorem} 

\begin{corollary}\label{cor:Robin-Stokes-tilde}
Under the assumptions of Theorem~\ref{thm:Robin-Stokes-Lip-tilde} and for $\del_0>$ large enough, the operator $\del_0+A_{B,q}$ has bounded imaginary powers. In particular, for $\alpha\in(0,1)$, we have
$$
 D((\del_0+{A}_{B,q})^\alpha)=[{L}^q_\si(\Om),D({A}_{B,q})]_\alpha.
$$
Moreover, the operator ${A}_{B,q}$ has maximal $L^p$-regularity, $p\in(1,\infty)$, on finite intervals in ${L}^q_\si(\Om)$. 
\end{corollary}

\begin{remark}\rm
(a) Again, Theorem~\ref{thm:Robin-Stokes-Lip-tilde} and Corollary~\ref{cor:Robin-Stokes-tilde} cover Stokes operators with Navier boundary conditions as in \eqref{eq:BC-alpha-beta} if we take $B$ as specified in \eqref{eq:Navier-Robin-BC}.

(b) The result on $L^p$-maximal regularity in Corollary~\ref{cor:Robin-Stokes-tilde} has been shown for Navier type boundary conditions in \cite{FarRost-MR}, but under an additional assumption on the uniform $C^{2,1}$-domain $\Om$.
\end{remark}

\begin{appendix}
\section{Auxiliary results}

\subsection{Traces and Gau{ss}'s theorem on unbounded domains}
We refer to \cite[Appendix B]{hobus-saal} for proofs of the following extensions of facts that are well-known for bounded domains. First we define, for any domain $\Om\subseteq\R^d$ and $q\in(1,\infty)$,
$$
 E_q(\Om):=\{f\in L^q(\Om)^d: \dv f\in L^q(\Om)\,\},
$$
which is a Banach space for $\nn{f}{E_q(\Om)}:=\nn{f}{L^q(\Om)^d}+\nn{\dv f}{L^q(\Om)}$.
If $\Om$ satisfies the segment property (so in particular if $\Om$ is a Lipschitz domain) then $C^\infty_c(\ov{\Om})^d$ is dense in $E_q(\Om)$ (see \cite[Lemma 13.1]{hobus-saal}). In the following proposition we collect the statements that are relevant for us.

\begin{prop}\label{prop:app-trace}
Let $\Om\subseteq\R^d$ be a uniform $C^{2,1}$-domain.
\begin{enumerate}[label=$(\roman*)$]
\item
For $q\in[1,\infty)$ the map $z\mapsto u|_{\partial\Om}$, defined on $C^\infty_c(\ov{\Om})$, has a continuous extension
$$
 \Tr:W^{1,q}(\Om) \to W^{1-\frac1q,q}(\partial\Om).
$$
For $q\in(1,\infty)$, $\Tr$ is surjective with a continuous linear right inverse 
$$\RTr:W^{1-\frac1q,q}(\partial\Om)\to W^{1,q}(\Om).$$
\item 
For any $u\in W^{1,1}(\Om)^d$ one has
$$
 \int_\Om \dv u\,dx = \int_{\partial\Om} \nu\cdot u\,d\si.
$$
\item 
For $q\in(1,\infty)$, $u\in W^{1,q}(\Om)$, and $v\in W^{1,q'}(\Om)^d$ one has
$$
  \int_\Om u\,\dv v\,dx=-\int_\Om\nabla u\cdot v\,dx + \int_{\partial\Om} u (\nu\cdot v)\,d\si.
$$
\item 
For $q\in(1,\infty)$ the map $v\mapsto \nu\cdot v|_{\partial\Om}$, defined on $C^\infty_c(\ov{\Om})$, has a continuous extension
$$
 \Tr_\nu: E_{q'}(\Om)\to W^{-\frac1{q'},q'}(\partial\Om)
 :=\big(W^{\frac1{q'},q}(\partial\Om)\big)'=\big(W^{1-\frac1{q},q}(\partial\Om)\big)',
$$
given by 
$$
 \sk{\Tr\, u}{\Tr_\nu v}_{\partial\Om}= \int_\Om u\,\dv v\,dx+\int_\Om\nabla u\cdot v\,dx
 \quad\mbox{for $u\in W^{1,q}(\Om)$.}
$$
Observe that $\sk{\Tr\,u}{\Tr_\nu v}_{\partial\Om}$ does not depend on the special choice of $u$ and we can take $u=\RTr\Tr\,u$. For simplicity of notation we put 
$$
 \sk{u}{\nu\cdot v}_{\partial\Om}:=\sk{\Tr\, u}{\Tr_\nu v}_{\partial\Om}\quad
 \mbox{for $u\in W^{1,q}(\Om)$ and $v\in E_{q'}(\Om)$.}
$$
\end{enumerate}
\end{prop}

For the proofs we refer to \cite[Lemmas B.2--B.7]{hobus-saal}. They may be extended to uniform Lipschitz domains.

\subsection{Extension, Sobolev embedding, and interpolation}
For the following extension operator we refer to \cite[Thm. VI.3.1/5]{Stein}. The formulation is the one from \cite[Lemma 12.2]{hobus-saal}.

\begin{prop}\label{prop:extension}
Let $\Om\subseteq\R^d$ be a uniform Lipschitz domain. Then there exists a linear operator $E$ mapping real-valued functions onto real-valued functions on $\R^d$ such that $Ef|_\Om=f$ holds for any function $f$ on $\Om$ and such that
$$
 E : W^{k,q}(\Om)\to W^{k,q}(\R^d)
$$
is bounded for all $1\le q<\infty$ and $k\in\N_0$.
\end{prop}

Using this extension operator $E$ one can prove the following Sobolev embeddings for $\Om$ via those on $\R^d$.

\begin{prop}\label{prop:sob-emb}
Let $\Om\subseteq\R^d$ be a uniform Lipschitz domain and $q\in(1,\infty)$ and $k\in\N$. 
If $q<\frac{d}{k}$ then $W^{k,q}(\Om)\emb L^r(\Om)$ where $\frac1r=\frac1q-\frac{k}{d}$.
If $q>\frac{d}k$ then $W^{k,q}(\Om)\emb C_0(\ov{\Om})$. 
\end{prop}

Using the extension operator $E$, the restriction $Rf=f|_\Om$, and \cite[1.2.4]{Triebel} one can also prove the following on complex interpolation spaces.

\begin{prop}\label{prop:app-interpol}
Let $\Om\subseteq\R^d$ be a uniform Lipschitz domain and $q\in(1,\infty)$. Then, for $k\in\N$ and $\theta\in(0,1)$,
$$
 [L^q(\Om),W^{k,q}(\Om)]_{\theta}=H^{\theta k,q}(\Om),
$$
where $H^{k\theta,q}(\Om)=R(H^{k\theta,q}(\R^d))$, i.e. the restrictions of functions in the Bessel potential space $H^{k\theta,q}(\R^d)$. If $k\theta=l\in\N$ then $H^{k\theta,q}(\Om)=W^{l,q}(\Om)$.
\end{prop}

As an application of Seeley's results (\cite{Seeley}) we obtain the following.

\begin{prop}\label{prop:seeley-interpol}
Let $\Om\subseteq\R^d$ be a uniform $C^{2,1}$-domain and $q\in(1,\infty)$. Then we have for 
$D(\Delta_{H,q})=\{u\in W^{2,q}(\Om)^d: \nu\cdot u=0, D_-(u)\nu=0\ \mbox{on $\partial\Om$}\}$ the following identities for complex interpolation spaces:
\begin{align*}
 [L^q(\Om)^d, D(\Delta_{H,q})]_\theta
 =\left\{\begin{array}{ll}
 H^{2\theta,q}(\Om)^d &, \ \theta\in(0,\frac1{2q}), \\
 \{u\in H^{2\theta,q}(\Om)^d: \nu\cdot u=0\ \mbox{on $\partial\Om$} \} &, \ \theta\in (\frac1{2q},\frac12+\frac1{2q}), \\
 \{u\in H^{2\theta,q}(\Om)^d: \nu\cdot u=0, D_-(u)\nu=0 \ \mbox{on $\partial\Om$} \} &, \ \theta\in (\frac12+\frac1{2q},1).
 \end{array}\right.
\end{align*}
\end{prop}

\begin{proof}
In order to apply the main result of \cite{Seeley} we rewrite the boundary condition $D_-(u)\nu=0$ in terms of normal derivatives of the components of $u$. Using \eqref{eq:weingarten} we obtain under the condition $\nu\cdot u=0$ on $\partial\Om$ that $D_-(u)\nu=0$ is equivalent to 
\begin{align*}
 (I-\nu\nu^T)\left[(\nabla u)^T\nu\right]=\left[(\nabla u)^T\nu\right]_{\tan}=\mathcal{W}u.
\end{align*}  
Hence we have exactly the form with the projection mentioned on p.54 before (3.4) in \cite{Seeley}. We can localize $\Om$ and apply then \cite[Theorem 4.1]{Seeley} using uniformity of $\Om$.
\end{proof}

\subsection{Generators in invariant subspaces}

The following lemma is easy. We include it with a proof for convenience of the reader.

\begin{lemma}\label{lem:inv-gen}
Let $X$ be a Banach space and $(T(t))_{t\ge0}$ be a $C_0$-semigroup in $X$ with negative generator~$A$. Let $Y$ be a closed subspace of $X$ that is invariant under each operator $T(t)$, $t\ge0$. 
Then $(S(t))_{t\ge0}:=(T(t)|_Y)_{t\ge0}$ is a $C_0$-semigroup in~$Y$
with negative generator $B=A|_{D(B)}$ where $D(B)=D(A)\cap Y$.
\end{lemma}

\begin{proof}
Clearly, $(S(t))_{t\ge0}$ is a $C_0$-semigroup in $Y$. If $y\in Y$ and $\frac1t(y-S(t)y)\to z$ in $Y$ then $\frac1t(y-T(t)y)\to z$ in $X$, and we conclude that $B$ is a restriction of $A$ and $D(B)\subseteq D(A)\cap Y$. If, on the other hand, $y\in Y$ and 
$\frac1t(y-T(t)y)\to z$ in $X$ then $z\in Y$ by closedness of $Y$, hence $\frac1t(y-S(t)y)\to z$ in $Y$ and $y\in D(B)$, $By=z$.
\end{proof}

We have the following corollary for fractional domain spaces.

\begin{corollary}\label{cor:inv-fract}
In the situation of Lemma~\ref{lem:inv-gen} let $\del\in\R$ be such that $(e^{-\del t}T(t))_{t\ge0}$ is bounded. For $\al\in(0,1)$ we then have $(\del+B)^\al=(\del+A)^\al|_{D((\del+B)^\al)}$ where $D((\del+B)^\al)=D((\del+A)^\al)\cap Y$.
\end{corollary}

\begin{proof}
Notice first that $\del+A$ is sectorial of angle $\le\frac\pi2$. Hence the fractional powers $(\del+A)^\al$ are well-defined and sectroial of angle $\le\al\frac\pi2$. In particular, $-(\del+A)^\al$ is the generator of a bounded analytic semigroup $(S_\al(t))$ and the semigroup operators may be represented by the holomorphic functional calculus of $A$ in terms of the resolvent operators of $A$. Since $Y$ is invariant under $(T(t))$, it is also invariant under the resolvents $(\la+A)^{-1}$ for $\mbox{Re}\la>\del$. We conclude that $Y$ is invariant under the semigroup $(S_\al(t))$. Then the assertion follows via Lemma~\ref{lem:inv-gen}.
\end{proof}

\end{appendix}


\small\normalsize\small

\begin{thebibliography}{99}

\bibitem{adams}
R.A. Adams, J.F.F. Fournier, Sobolev spaces, Second edition, Pure and Applied Mathematics (Amsterdam), 140. Elsevier/Academic Press, Amsterdam, 2003.

\bibitem{AHS}
H. Amann, M. Hieber, G. Simonett, Bounded $H_\infty$ calculus for elliptic
operators, Differential Integral Equations 7, 613--653 (1994).

\bibitem{Abe:cyl}
K. Abe,
The Navier–Stokes equations with the Neumann boundary condition in an infinite cylinder
manuscripta math. 160, 359--383 (2019).

\bibitem{AlBaba-NBC}
Hind Al Baba,
Maximal $L^p$-$L^q$ regularity to the Stokes
problem with Navier boundary conditions,
Adv. Nonlinear Anal. 8, 743--761 (2019).

\bibitem{AKST}
T. Akiyama, H. Kasai, Y. Shibata, M. Tsutsumi, On a resolvent estimate of a system of Laplace operators with perfect wall condition,  
Funkcial. Ekvac. 47, no. 3, 361--394 (2004). 

\bibitem{BabK}
J. Babutzka, P.C. Kunstmann, 
$L^q$-Helmholtz decomposition on periodic domains and applications to Navier-Stokes equations,
J. Math. Fluid Mech. 20, no. 3, 1093-1121 (2018).

\bibitem{BL:interpol}
J. Bergh, J. L\"ofstr\"om, 
Interpolation Spaces. An Introduction, Grundlehren der mathematischen 
Wissenschaften 223, Berlin-Heidelberg-New York, Springer-Verlag, 1976.

\bibitem{BK-mr}
S. Blunck, P.C. Kunstmann, 
Weighted norm estimates and maximal regularity, 
Adv. Differential Equations 7, no. 12, 1513--1532 (2002).

\bibitem{BK-hinfty}
S. Blunck, P.C. Kunstmann,
Calder\'on-Zygmund theory for non-integral operators and the $H^\infty$-functional calculus. Rev. Mat. Iberoamericana 19, no. 3, 919–942 (2003).

\bibitem{Bog86}
M. E. Bogovski\u{\i}, Decomposition of $L_p(­\Omega;\R^n)$ into a direct sum of subspaces of solenoidal
and potential vector fields, Dokl. Akad. Nauk SSSR, 286 (1986), 781--786.

\bibitem{calderon}
A.P. Calder\'on, 
Intermediate spaces and interpolation, the complex method,
Stud. Math. 24, 113--190 (1964).

\bibitem{DOS}
X.T. Duong, E.M. Ouhabaz, A. Sikora, 
Plancherel-type estimates and sharp spectral multipliers. J. Funct. Anal. 196 , no. 2, 443--485 (2002).

\bibitem{duong-robinson}
X.T. Duong, D.W. Robinson, 
Semigroup kernels, Poisson bounds, and holomorphic functional calculus, 
J. Funct. Anal. 142, no. 1, 89--128 (1996).

\bibitem{FKS:ArchMath}
R. Farwig, H. Kozono, H. Sohr,
On the Helmholtz decomposition in general unbounded domains,
Arch. Math. 88, 239--248 (2007).

\bibitem{FarRost-Res}
R. Farwig, V. Rosteck,
Resolvent estimates of the Stokes system with Navier boundary conditions in general unbounded domains, 
Adv. Differential Equations 21, no. 5-6, 401--428 (2016).

\bibitem{FarRost-MR}
R. Farwig, V. Rosteck, Maximal regularity of the Stokes system with Navier boundary condition in general unbounded domains, 
J. Math. Soc. Japan 71, no. 4, 1293--1319 (2019).

\bibitem{GHHS}
M. Geissert, H. Heck, M. Hieber, O. Sawada, Weak Neumann implies Stokes, 
J. Reine Angew. Math. 669, 75--100  (2012).

\bibitem{GHT}
M. Geissert, H. Heck, C. Trunk, $H^\infty$-calculus for a system of Laplace operators with mixed order boundary conditions, Discrete Contin. Dyn. Syst. Ser. S 6, no. 5, 1259--1275 (2013).

\bibitem{GK}
M. Geissert, P.C. Kunstmann, Weak Neumann implies $H^\infty$ for Stokes,
J. Math. Soc. Japan 67, no. 1, 183--193 (2015).

\bibitem{Giga}
Y. Giga, Domains of fractional powers of the Stokes operator in $L_r$ spaces,
Arch. Ration. Mech. Anal. 89, 251--265 (1985).

\bibitem{Greiner}
G. Greiner, Perturbing the boundary conditions of a generator,
Houston J. Math. 13, no. 2, 213--229 (1987).

\bibitem{HK:nse}
B.H. Haak, P.C. Kunstmann, On Kato's method for Navier Stokes equations,  
J. Math. Fluid Mech. 11, no. 4, 492--535 (2009).

\bibitem{fc-book}
M. Haase, The Functional Calculus for Sectorial Operators, Operator Theory: Advances and Applications Vol. 169, Birkh\"auser, 2006.

\bibitem{hobus-saal}
P. Hobus, J. Saal, Stokes and Navier-Stokes equations subject to partial slip on uniform $C^{2,1}$-domains in $L^q$-spaces, 
J. Differential Equations 284, 374--432 (2021).

\bibitem{KKW}
N.J. Kalton, P.C. Kunstmann, L. Weis,
Perturbation and interpolation theorems for the $H^\infty$-calculus with
applications to differential operators,
Math. Ann. 336, 747--801 (2006).

\bibitem{KM}
N.J. Kalton, M. Mitrea, 
Stability results on interpolation scales of quasi-Banach spaces and applications, 
Trans. Am. Math. Soc. 350, 3903-3922 (1998).

\bibitem{K:london}
P.C. Kunstmann, 
Heat kernel estimates and $L^p$ spectral independence of elliptic operators, 
Bull. London Math. Soc. 31, no. 3, 345--353 (1999).

\bibitem{pck-m2bc}
P.C. Kunstmann, 
Maximal $L^p$-regularity for second order elliptic operators with
uniformly continuous coefficients on domains, 
in Iannelli, Mimmo (ed.) et al., Evolution equations: applications to 
physics, industry, life sciences and economics, 
Basel, Birkh\"auser. Prog. Nonlinear Differ. Equ. Appl. 55, 293-305 (2003).

\bibitem{K:hinfty-stokes}
P.C. Kunstmann, $H^\infty$-calculus for the Stokes operator on unbounded domains, 
Arch. Math. 91, no. 2, 178--186 (2008).

\bibitem{KUhl:spec-mult}
P.C. Kunstmann, M. Uhl,
Spectral multiplier theorems of H\"ormander type on Hardy and Lebesgue spaces, 
J. Operator Theory 73, no. 1, 27--69 (2015).

\bibitem{KUhl:ell-sys}
P.C. Kunstmann, M. Uhl, 
$L^p$-spectral multipliers for some elliptic systems,
Proc. Edinb. Math. Soc. (2) 58, no. 1, 231--253 (2015).

\bibitem{KuW:Pisa}
P.C. Kunstmann, L. Weis, 
Perturbation theorems for maximal $L_p$-regularity, 
Ann. Scuola Norm. Sup. Pisa Cl. Sci. (4) 30, no. 2, 415--435 (2001).

\bibitem{KuW:levico}
P.C. Kunstmann, L. Weis,
Maximal $L_p$-regularity for parabolic equations, Fourier multiplier theorems and $H^\infty$-functional calculus,
in M. Iannelli, R. Nagel, S. Piazzera (eds.), \emph{Functional Analytic Methods for Evolution Equations},
Springer Lecture Notes Math. 1855, 65--311 (2004).

\bibitem{KuW:erratum}
P.C. Kunstmann, L. Weis, 
Erratum to: Perturbation and interpolation theorems for the $H^\infty$-calculus with applications to differential operators, 
Math. Ann. 357, no. 2, 801--804 (2013).

\bibitem{KuW:Hinfty-Stokes}
P.C. Kunstmann, L. Weis, 
New criteria for the $H^\infty$-calculus and the Stokes operator on bounded Lipschitz domains, 
J. Evol. Equ. 17, no. 1, 387--409 (2017).

\bibitem{mas-bog-appr}
V.N. Maslennikova, M.E. Bogovskii, Approximation of potential and solenoidal vector fields, Sibirsk. Mat. Zh. 24, no. 5, 149--171 (1983). 

\bibitem{MaBo}
V.N. Maslennikova, M.E. Bogovskii, Elliptic boundary value problems in
unbounded domains with noncompact and nonsmooth boundaries. Rend. Sem. Mat.
Fis. Milano. LVI, 125--138 (1986).

\bibitem{MM:HNS-nonlin}
M. Mitrea, S. Monniaux. The nonlinear Hodge-Navier-Stokes equations in Lipschitz domains, 
Differential and Integral equations 22 (3-4), 339--356 (2009).

\bibitem{MM:TAMS}
M. Mitrea, S. Monniaux, On the analyticity of the semigroup generated by the Stokes operator with Neumann-type boundary conditions on Lipschitz subdomains of Riemannian manifolds, Trans. Amer. Math. Soc. 361 (no. 6), 3125--3157 (2009).

\bibitem{Mon-Ou}
S. Monniaux, E.M. Ouhabaz, 
The incompressible Navier-Stokes system with time-dependent Robin-type boundary conditions,
J. Math. Fluid Mech. 17, no. 4, 707--722 (2015).

\bibitem{NoSa}
A. Noll, J. Saal, $H^\infty$-calculus for the Stokes operator on $L^q$-spaces,
Math. Z. 244, no. 3, 651--688 (2003).

\bibitem{Ouh:inv-convex}
E.M. Ouhabaz, Invariance of closed convex sets and domination criteria for semigroups,
Potential Analysis 5, 611--625 (1996). 

\bibitem{Ouh:book}
E.M. Ouhabaz, 
Analysis of heat equations on domains, London Mathematical Society Monographs Series 31, Princeton University Press, Princeton, NJ, 2005.

\bibitem{Pruess}
J. Pr\"uss, $H^\infty$-calculus for generalized Stokes operators,
J. Evol. Equ. 18, no. 3, 1543--1574 (2018).

\bibitem{Seeley}
R. Seeley, Interpolation in $L^p$ with boundary conditions, 
Studia Math. 44, 47--60 (1972). 

\bibitem{Sohr}
H. Sohr, The Navier-Stokes equations. Birkh\"auser/Springer Basel AG, Basel, 2001.

\bibitem{solonnikov}
V.A. Solonnikov, On the solvability of boundary and initial-boundary value problems for the Navier-Stokes system in domains with noncompact boundaries, 
Pacific J. Math. 93, no. 2, 443--458 (1981).

\bibitem{Stein}
E. M. Stein. Singular integrals and differentiability properties of functions. Princeton University Press,
Princeton, N.J., 1970.

\bibitem{Triebel}
H. Triebel, Interpolation, Function Spaces, Differential Operators, 
North-Holland Mathematical Library. Vol. 18. Amsterdam - New York - Oxford,
1978.

\end{thebibliography}
\end{document}